\documentclass[11pt]{article}
\usepackage{pgfplots}
\usepackage{amsmath}
\usepackage{amssymb}
\usepackage{algorithmic}
\usepackage{algorithm}
\usepackage{mathtools}
\usepackage{color}
\usepackage{tikz}
\usetikzlibrary{arrows,
	arrows.meta,
	bending,calc,patterns,angles,quotes}
\usepackage{float}

\newcommand{\rd}{{\rm d}}

\makeatletter
\providecommand*{\cupdot}{%
	\mathbin{%
		\mathpalette\@cupdot{}%
	}%
}
\newcommand*{\@cupdot}[2]{%
	\ooalign{%
		$\m@th#1\cup$\cr
		\hidewidth$\m@th#1\cdot$\hidewidth
	}%
}
\makeatother

\newtheorem{theorem}{Theorem}

\newtheorem{lemma}{Lemma}
\newtheorem{proposition}{Proposition}
\newtheorem{remark}{Remark}

\newtheorem{corollary}{Corollary}
\def\endpf{{\ \hfill\hbox{\vrule width1.0ex height1.0ex}\parfillskip 0pt
	}}
	
\begin{document}
	\title{Externalities in queues as stochastic processes: \newline The case of FCFS M/G/1}
	\author{Royi Jacobovic and Michel Mandjes \thanks{Korteweg-de Vries Institute; University of Amsterdam; 1098 XG Amsterdam; Netherlands. 
			{\tt royi.jacobovic@mail.huji.ac.il, M.R.H.Mandjes@uva.nl}. This research was supported by the European Union’s Horizon 2020 research and innovation programme under the Marie Skłodowska-Curie grant agreement no. 945045, and by the NWO Gravitation project NETWORKS under grant no. 024.002.003. \includegraphics[height=1em]{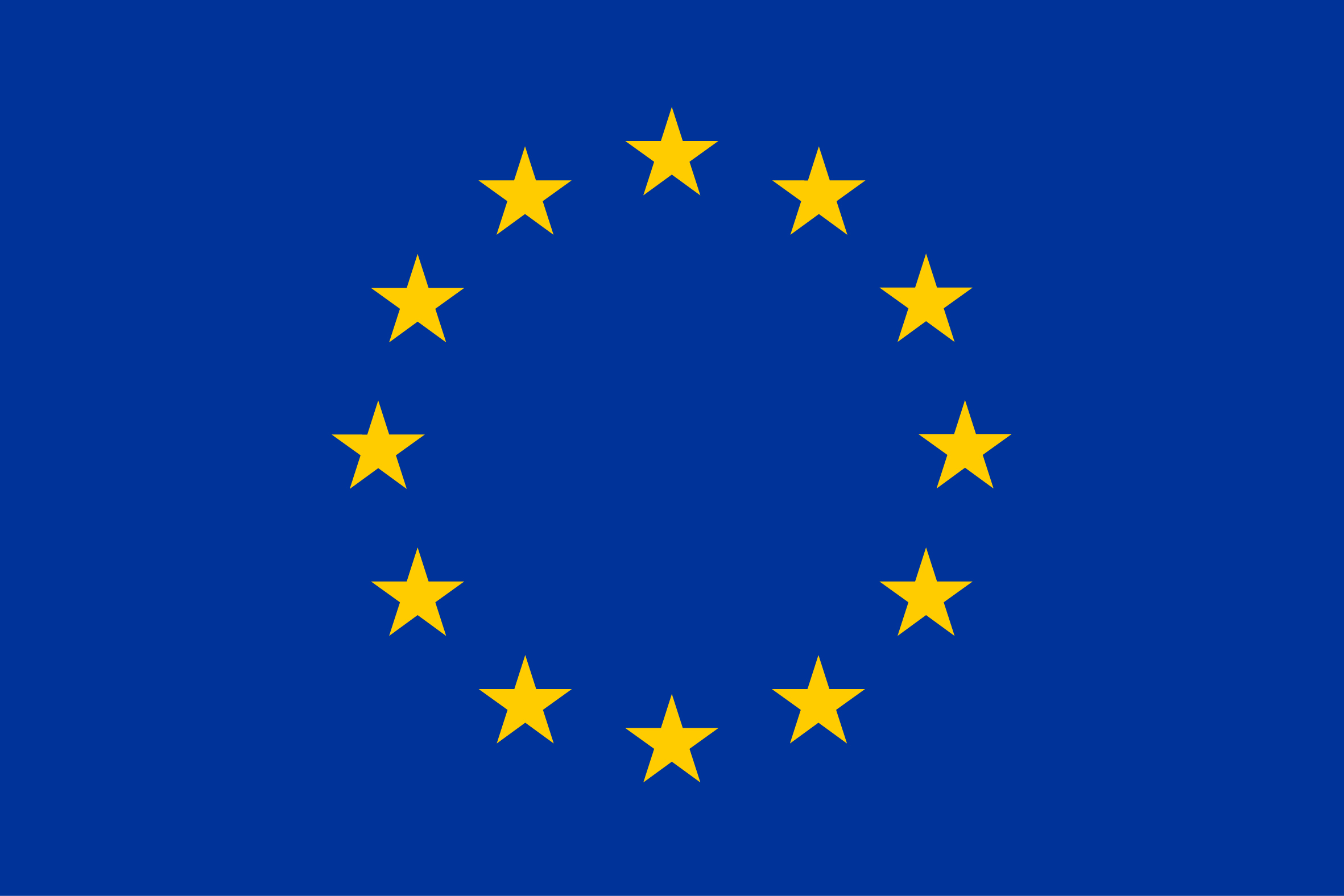} } }
	\date{\today}
	\maketitle
	\begin{abstract} \noindent
		Externalities are the costs that a user of a common resource imposes on others. In the context of an FCFS M/G/1  queue, where a customer with service demand $x\geq0$ arrives  when the workload level is $v\geq0$, the externality $E_v(x)$ is the total waiting time which could be saved if this customer gave up on their service demand. In this work, we analyze the \textit{externalities process} $E_v(\cdot)=\left\{E_v(x):x\geq0\right\}$. It is shown that this process can be represented by an integral of a (shifted in time by $v$) compound Poisson process with a positive discrete jump distribution, so that $E_v(\cdot)$ is convex. Furthermore, we compute the Laplace-Stieltjes transform (LST) of the finite-dimensional distributions of $E_v(\cdot)$ as well as its mean and auto-covariance functions. We also identify conditions under which a sequence of normalized externalities processes admits a weak convergence on $\mathcal{D}[0,\infty)$ equipped with the uniform metric to an integral of a (shifted in time by $v$) standard Wiener process. Finally, we also consider the extended framework when $v$ is a general nonnegative random variable which is independent from the arrival process and the service demands. Our analysis leads to substantial generalizations of the results presented in the seminal paper by Haviv and Ritov (1998).      
	\end{abstract}
	
	
	\bigskip
	\noindent {\bf AMS Subject Classification (MSC2010):}  60K25, 60K30, 60K37.
	\section{Introduction}\label{sec: introduction}
	Consider a conventional M/G/1 queueing system that is served according to the first-come, first-served (FCFS) discipline, with arrival rate $\lambda>0$ and with the service distribution given by $B(\cdot)$. Assume that the queue is stable, and let the workload level at time $t=0$ be $v\geq0$ (say) minutes. Denote the workload at time $t\geq0$ by $W_v(t)$, and let $T_i$ be the arrival time of the $i$-th customer. The main objective of this paper is to analyze the aggregate effect of an additional customer, who has arrived at time $t=0$ with a service requirement of size $x\geq0$, on the waiting times of {\it all} other customers. In other words, we are interested in the distribution of the {\it externality}:
	\begin{equation}\label{eq: externalities definition1}
	E(x,v)\equiv\sum_{i=1}^\infty\left[W_{v+x}(T_i)-W_v(T_i)\right].
	\end{equation}
	Thus, the externality $E(x,v)$ is to be interpreted as the the total waiting time which could be saved if the additional customer reduced their service requirement from $x\geqslant 0$ to zero. To the best of our knowledge, \cite{Haviv1998} is the only existing paper that analyzes $E(x,v)$. In \cite{Haviv1998} it was shown that if (i)~$v$ is a random variable which is independent from the arrival process and the service requirements of the customers, and (ii) $v$ is distributed according to the stationary distribution of the workload process, then the mean of $E(x,v)$ is given by
	\begin{equation}\label{eq: haviv}
	    \mathbb{E}\left[E(x,v)\right]=\frac{\lambda x}{2\left(1-\rho\right)}\left(\frac{\lambda\mu_2}{1-\rho}+x\right)
	\end{equation}
	where $\mu_i$ ($i=1,2,\ldots$) is the $i$-th moment pertaining to $B(\cdot)$ and $\rho\equiv\lambda\mu_1$.

 \medskip

 Whereas \cite{Haviv1998} focused on computing the mean externality (under the specific condition mentioned above), we have managed to develop a full probabilistic analysis of $E(x,v)$. In this context it is important to notice that 
		\begin{equation}
		    \{E(x,v): x,v\geq0\}
		\end{equation} can be seen as a collection of random variables which are all defined on the same probability space. By considering $v$ as a fixed parameter while $x$ is given the role of a time index, we analyze the stochastic process $E_v(\cdot)\equiv E(\cdot,v)$. To underline the natural interpretation of this process, let the additional customer arrive to the queue when the existing workload is $v$, and assume that this customer has {\it two} tasks that they want the server to do for them: a first one of size $x_1\geq 0$ and a second one of size $x_2\geq 0$. Then, $E_v(x_1+x_2)-E_v(x_1)$ is equal to the total waiting time that could be saved by the other customers if the customer gave up on their second task but insisted on completing the first one.

  \medskip
	
	The main contribution of this work lies in an extensive analysis of $E_v(\cdot)$ which in the sequel we refer to as the \textit{the externalities process}. Specific open questions which we have managed to solve in the current paper, are: 
		\begin{enumerate}
			\item What can be said about the distribution of the externalities in a {\it non-stationary} FCFS M/G/1 queue? As it turns out, the externalities process $E_v(\cdot)$ can be represented by an integral of a compound Poisson process that is shifted in time by an amount $v$. Importantly, this compound Poisson process is defined on the same probability space as the one on which our model is defined.
			
			\item Observe that the expected value in \eqref{eq: haviv} is convex in $x$, which indicates that the marginal effect of extra workload on the customer population is increasing. Is it possible to extend this result by showing convexity of the  externalities process? The answer is affirmative, where we also provide an explicit representation of the corresponding right-derivative.
			
			\item Is there a systematic way to evaluate the moments of the externality $E_v(x)$? To this end, we derive the Laplace-Stieltjes transform (LST) of the finite-dimensional distributions of the externalities process $E_v(\cdot)$, from which the moments follow. In particular, we provide closed-form formulae for the auto-covariance and auto-correlation functions of $E_v(\cdot)$. Remarkably, it is shown that when $v$ is fixed, then the auto-correlation does not depend on the stochastic ingredients of the model, i.e., the arrival rate and service distribution. 
			
			\item Is it possible to approximate the distribution of the externalities in some asymptotic regime? We show that, under an appropriate scaling, there is convergence of $E_v(\cdot)$ to a specific Gaussian limiting process. The convergence takes place as the arrival rate tends to infinity and the service distribution is `well behaved', e.g., it tends to zero in an appropriate way. 
		\end{enumerate}     
		
	\subsection{Motivation}
 We proceed by discussing the relevance  of our result, and their applicability in an operational context. We do so by distinguishing three strands of application domains. 

 \medskip

 \noindent
	{\bf Choice of a management scheme.}
	 In the introduction of their paper, Haviv and Ritov~\cite{Haviv1998} discuss various applications of the externalities setup that they analyze: airplanes taking off from a runway, commuters crossing a bridge, jobs sharing a common CPU, and  messages being routed through a common data network. Their motivation for studying the distribution of externalities is as follows. In the first place, they argue that ``a zero profit operator who charges users for the use of a common facility usually likes to do so in accordance with the congestion costs that they impose on others". This aligns with results in e.g.\ \cite{Ha2001,Haviv2014,Haviv2018a,Haviv2018b,Jacobovic2022} where 
	 various relations between optimal queue regulation schemes and externalities are revealed. Then, they point out that there are various policies of managing a queueing system (e.g., by implementing different service disciplines). Correspondingly, different management policies may result in different amounts of externalities imposed by the same user. This leads them to the conclusion that ``the resulting pricing mechanism can serve as an additional criterion for deciding which management scheme to adopt". 
  A general account of externalities in a queueing context is given in \cite{HH}, as well as various other connections between queueing and game theory.

  \medskip

 \noindent
	{\bf Queues with discretionary services.}\label{subsubsec: discretionary services}
	Recently, there has been a growing interest in queueing models with customers who {\it themselves} choose their service durations (see, e.g.,\ \cite{Feldman2022,Jacobovic2022} and the references therein). When considering single-server queues with a non-preemptive service discipline, the customer who gets service does not care about the increasing costs of the waiting customers behind them, thus yielding a resource allocation which is inefficient from a social point of view. In order to restore social efficiency, a social planner may want to impose some sort of regulation. For example, the planner may decide on a price function which tells every customer how much they are going to pay for every service duration to be purchased. A price function will be optimal if it makes the customers behave as they should according to the socially optimal resource allocation.  A reasonable price mechanism amounts to requiring every customer to pay for the expected cost which is enforced on the others due to their service requirement.  
	
	The earlier paper \cite{Jacobovic2022}  considered a model of a single-server queue with customers who arrive according to a Poisson process and dynamically choose their service durations, showing that when the social planner is restricted to choose a price function which is determined by the service requirement only, then the optimal price function internalizes the (expected) externalities. It is an open problem  \cite{Jacobovic2022a} whether a similar phenomenon occurs when the social planner may choose a price function which depends on the state of the queue. If the answer to this question is affirmative and the social planner observes the workload level at the onset of every service duration, then the optimal price function is equal to $\mathbb{E}[E_v(x)]$, with $v$ the initial workload {at the start of the service} and $x$ the corresponding  service requirement. As is shown in the present paper, this would reduce the search for the optimal price function to the parametric family of quadratic functions in $x$ which are also linear in $v$. 
	
	Similarly, in another possible scenario a social planner observes the number of waiting customers at the start of every service but they do not see the customers' service requirements (See also \cite[Section 3]{Haviv2014}). In this case, the conjectured optimal price function is the conditional expectation of $E_v(x)$ given the available information at the start of the service. Once more, our results imply that this conditional expectation is quadratic in $x$ and also linear in the number of waiting customers at the time of the start of the service. 

 \medskip
	
	\noindent
	{\bf Queues with a proactive service discipline.}
	Consider an emergency room with a single specific bed which is reserved for patients with special needs, e.g., those who arrive because of strokes, heart attacks, etc. We refer to these patients as `urgent', while the patients who arrive due to other reasons are called `regular'. Note that the special bed might be useful also for regular patients while the urgent ones can be treated only in the special bed. Hence a non-trivial question is: if there are many regular patients and no urgent patients, should the regular patients be allowed to use the special bed? Doing so is evidently beneficial to the regular patients, but it is also possible that immediately after allocating a regular patient to the special bed, a batch of urgent patients arrives whose treatments will be delayed. 

    Now, assume that the urgent patients arrive according to a Poisson process with rate $\lambda$ and their service requirements are iid random variables with a distribution function $B(\cdot)$ which are independent from the arrival process. Then, observe that $E_0(x)$ is equal to the total damage which is caused to the urgent patients due to an allocation of a regular customer into the special bed for $x$ minutes once it is empty. Clearly, the decision maker could benefit from the distributional properties of $E_0(x)$ that we establish in the present paper.
    
    The above example connects our work with  server-allocation problems in multiclass queues. Recent progress in this direction can be found in, e.g., \cite{Chan2021,Hu2022,Huang2015,Liu2022}.
	
	\subsection{Organization of the paper}	
 The organization of this work is as follows. Section \ref{sec: busy period} starts by a brief discussion of a known result, extensively used in the papwer: a fixed-point relation which is satisfied by the LST of the distribution of the number of customers who arrive to a queue during a busy period. Besides this fixed-point relation, all results presented are novel contributions.
 Then, Section~\ref{sec: definition} includes a representation of the externalities process $E_v(\cdot)$ in terms of a  compound Poisson process,  yielding two insightful decompositions:

 \medskip
 
	\noindent{\bf Decomposition 1.} $E_v(\cdot)$ is equal to an integral of a compound Poisson process which is shifted in time by $v$.  The rate of this process is equal to $\lambda$ and its jumps have the distribution identified in Section~\ref{sec: busy period}.  Section~\ref{sec: crossing times} provides a compact analysis of the crossing times of the right-derivative of $E_v(\cdot)$. An important application of this decomposition can be found in Section \ref{sec: heavy traffic} where we derive of a functional central limit theorem for the externalities process.

 \medskip
	
	\noindent{\bf Decomposition 2.} The distribution of $E_v(x_1+x_2)-E_v(x_1)$ is the same as the distribution of a sum of independent random variables. This helps in Section \ref{sec: moments} where we derive the LST of the finite-dimensional distributions pertaining to the process $E_v(\cdot)$. Moreover, this decomposition plays an important role in the derivations in Section \ref{sec: v is random} where we consider the more general framework when $v$ is a nonnegative random variable, independent from the arrival process and the service requirements of the customers. In particular, the results of this part include a generalization of \eqref{eq: haviv} to the case where $v$ is not necessarily distributed according to the stationary distribution of the workload process. 
    
    \medskip

    \noindent
    Section~\ref{sec: conclusion} concludes by discussing some related open problems which lead to several directions of future research. In order to optimize the flow of the paper,  all proofs are given in Section~\ref{sec: proofs}. 
	
	\section{Number of customers during busy period}\label{sec: busy period}
	This section discusses a few results concerning the number of customers who arrive to a stable FCFS M/G/1  queue during a single busy period, needed  in the upcoming sections. Proposition \ref{prop: recursive relation} is standard \cite[Chapter II.4.4]{Cohen1969}, while all the other results in this section are essentially direct consequences. However, since we did not find a reference for Propositions~\ref{prop: recursion}--\ref{prop: PGF}, we decided to include their proofs. For additional work on the distribution of the number of customers who arrive during a busy period, see \cite{Novak2006} and the references therein. 
	
	As before, we consider the setting of an M/G/1 queue with arrival rate $\lambda$ and a service distribution $B(\cdot)$, but now the system starts empty at time $t=0$. In addition, denote the LST of $B(\cdot)$ by \[b(s)\equiv \int_0^\infty e^{-st}\rd B(t),\:\:\:\:\:s>0.\] 
	and, for any $n\geq 1$, denote the $n$-th moment of $B(\cdot)$ by	
	\begin{equation}
	\mu_n\equiv\int_0^\infty t^n\rd B(t)\,.
	\end{equation} 
	Throughout this paper we assume that $\rho\equiv\lambda\mu_1\in(0,1)$ to ensure stability.

	Let $N(s)$ be the probability that exactly $s$ customers received service during the first busy period.
	The associated $k$-th moment is denoted by
	\begin{equation}
		\eta_k\equiv\sum_{s=1}^\infty s^kN(s)\ \ , \ \ k\geq1\,.
	\end{equation} 
	
	\begin{proposition}\label{prop: recursive relation}
		For every $z\in(0,1)$, the following fixed-point equation in $y$
		\begin{equation}\label{eq: LST fixed point}
			y=z\,b\big(\lambda(1-y)\big)\,.
		\end{equation}
		has a unique solution $y_z$ which belongs to $(0,1)$. Furthermore, $y_z$ equals the generating function  
		\begin{equation}
			\hat{N}(z)\equiv\sum_{s=1}^\infty z^s N(s)\,.
		\end{equation}
	\end{proposition}
	
	\begin{remark}
		\normalfont Notice that 
		\begin{equation}
			0<zb\left(\lambda\right)<zb(0)<1\,.
		\end{equation}
		Thus, since both sides of \eqref{eq: LST fixed point} are continuous in $y$, for every $z\in(0,1)$, it is possible to find $y_z$ efficiently by a standard line-search algorithm.
	\end{remark} 
	
	In particular, for every $\alpha>0$, we can insert $z=e^{-\alpha}$ into \eqref{eq: LST fixed point}. This yields the following fixed-point relation for the LST:
	\begin{equation}\label{eq: LST}
		\widetilde{N}(\alpha)\equiv\sum_{s=1}^\infty e^{-\alpha s}N(s)=e^{-\alpha}b\left\{\lambda\left[1-\widetilde{N}(\alpha)\right]\right\}\ \ , \ \ \alpha>0\,.
	\end{equation} 
	Therefore, we can differentiate both sides of \eqref{eq: LST} at zero in order to get a recursive formula for the moments $\eta_n$, $n\geq1$. In the sequel, for any pair of integers $m$ and $k$ such that $1\leq m\leq k$, denote the corresponding incomplete Bell's polynomial
	\begin{align}
		\mathcal{B}_{k,m}&[x_1,x_2,\ldots,x_{k-m+1}]\\&\equiv\sum\frac{k!}{j_1!j_2!\ldots j_{k-m+1}!}\left(\frac{x_1}{1!}\right)^{j_1}\left(\frac{x_2}{2!}\right)^{j_2}\ldots\left(\frac{x_{k-m+1}}{(k-m+1)!}\right)^{j_{k-m+1}}\nonumber
	\end{align} 
	where the summation is over all non-negative integers $j_1,j_2,\ldots,j_{k-m+1}$ which satisfy the following two conditions:
	\begin{equation}
		\sum_{i=1}^{k-m+1}j_i=m\ \ , \ \ \sum_{i=1}^{k-m+1}ij_i=k\,.
	\end{equation}
	In addition, for any pair of integers $m$ and $k$ such that $1\leq m\leq k$, we introduce the following compact notation:
	\[\check{\mathcal{B}}_{k,m}\equiv \mathcal{B}_{k,m}\left[-\eta_1,\eta_2,\ldots,(-1)^{(k-m+1)}\eta_{k-m+1}\right].\]
	\begin{proposition}\label{prop: recursion}
		For every positive integer $n$, 
		\begin{equation}\label{eq: moment}
			\eta_n\,=\frac{(-1)^n}{1-\rho}\left\{(-1)^n+\sum_{k=1}^{n-1}\binom{n}{k}(-1)^{n-k}\sum_{m=1}^k\lambda^m\mu_m\check{\mathcal{B}}_{k,m}+\sum_{m=2}^n\lambda^m\mu_m\check{\mathcal{B}}_{n,m}\right\}\,.
		\end{equation}
	\end{proposition}
	
	The following corollary, providing explicit expressions for the first three moments in terms of the moments of $B(\cdot)$, is an immediate consequence of Proposition \ref{prop: recursion}. The first moment
		$\eta_1$ also follows from the well-known result that the expected length of the busy period is $\mu_1/(1-\rho)$, in combination with Little's law.
	\begin{corollary}\label{cor: N moments} The first three moments are given by
		\begin{align}
			&\eta_1=\frac{1}{1-\rho}\,, \\&\eta_2=\frac{1}{1-\rho}\cdot\left[1+\frac{2\rho}{1-\rho}+\frac{\lambda^2\mu_2}{(1-\rho)^2}\right]\,,\nonumber\\&\eta_3=\frac{1}{1-\rho}\left\{1+\frac{3\rho}{1-\rho}+3\left[\rho\eta_2+\frac{\lambda^2\mu_2}{(1-\rho)^2}\right]+\frac{3\lambda^2\mu_2\eta_2}{1-\rho}+\frac{\lambda^3\mu_3}{(1-\rho)^3}\right\}.\nonumber
		\end{align}
	\end{corollary}

In a similar fashion, a combinatorial formula for the probability mass function $N(s)$, $s=1,2,\ldots$ may be derived by repeatedly differentiating
	\begin{equation}\label{eq: PGF}
		\hat{N}(z) =zb\left\{\lambda\left[1-\hat{N}(z)\right]\right\}\ \ , \ \ |z|<1
	\end{equation} 
	at zero. Using the compact notation
	\[\bar{\mathcal{B}}_{s-1,m}\equiv \mathcal{B}_{s-1,m}\left[N(1),2N(2),\ldots,(s-m)!N(s-m)\right],\]
	we arrive at the following recursion.
	\begin{proposition}\label{prop: PGF}
		$N(1)=b(\lambda)$ and for every $s\geq2$,
		
		\begin{equation}
			N(s)=\frac{1}{(s-1)!}\sum_{m=1}^{s-1}(-\lambda)^m\,b^{(m)}(\lambda)\,\bar{\mathcal{B}}_{s-1,m}\,.
		\end{equation}
	\end{proposition}
	
	\section{Decompositions of externalities}\label{sec: definition}
	
	This section first introduces the notation that will be used throughout the paper, and provides a detailed model description. Then we state our decomposition results.
	\subsection{Model description}
	With $\lambda$ and $B(\cdot)$ as defined before, let $\left\{J(t):t\geq0\right\}$ be a compound Poisson process with rate $\lambda\in(0,\infty)$ and a nonnegative jump distribution $B(\cdot)$. In addition, for each $i\geq1$, we let $T_i$ be the time of the $i$-th jump of the process $J(\cdot)$. In addition,  consider two processes $X_1(\cdot)$ and $X_2(\cdot)$ which are  given by
	\begin{equation}
	X_1(t)\equiv X_2(t)-x\equiv v+J(t)-t\ \ , \ \ t\geq0\,,
	\end{equation}
	for some two parameters $x,v\geq0$. Then, for each $i=1,2$, let $Y_i(\cdot)$ be the reflection of $X_i(\cdot)$ at the origin; this reflection, formally defined in e.g.\  \cite[Section 2.4]{Debicki2015}, can be thought of as a mechanism preventing the `free processes' $X_i(\cdot)$ from becoming negative. 
	Then, define, for a given initial workload $v$ and service requirement $x$, the externality via \begin{equation}\label{eq: externalities definition}
	E(x,v)\equiv\sum_{n=1}^\infty\left[Y_2(T_n)-Y_1(T_n)\right].
	\end{equation}
	Notice that $Y_2(t)\geq Y_1(t)$, but the stability condition
	$\rho<1$ implies that the hitting time of $Y_2(\cdot)$ in the origin is an almost surely finite random variable. Denote this random  variable by $\zeta$ and notice that this makes $E(x,v)$ an almost surely finite random variable. Observe that from time $\zeta$ on, the processes $Y_1(t)$ and $Y_2(t)$ are {\it coupled} (in that they coincide). 
	
	Importantly, $Y_1(\cdot)$ (resp.\ $Y_2(\cdot)$) coincides with $W_v(\cdot)$ (resp.\ $W_{v+x}(\cdot)$) which was defined in the beginning of Section \ref{sec: introduction}. Therefore, the quantity $E(x,v)$ represents the externality which is due to an arrival of a customer with a service demand of $x$ when the processing time of the existing workload is $v$. More generally, fixing the initial workload $v\geq0$, we can consider a stochastic process $E_v(x)\equiv E(x,v)$ indexed by $x\in[0,\infty)$, which in the sequel we refer to as the \textit{externalities process}. 
	
	\subsection{Decomposition 1}
	For the analysis of the externalities process, the following notation and definitions are needed. Throughout, the initial workload $v$ pertaining to $Y_1(t)$ is held fixed. In the first place, let $\tau_0$ be the end of the first busy period of $Y_1(t)$. Also, let $\sigma_1$ be the time of the first jump of $J(\cdot)$ which occurs after $\tau_0$. In addition, denote the first time after $\sigma_1$ in which $Y_1(\cdot)$ hits the origin by $\tau_1$ (i.e., the end of the second busy period of $Y_1(t)$). Similarly, we can define $\sigma_2$ to be the time of the first jump of $J(\cdot)$ which occurs after $\tau_1$. Moreover, let $\tau_2$ be the first time after $\sigma_2$ in which $Y_1(\cdot)$ hits the origin. We may continue recursively with this construction in the evident manner, thus yielding the two sequences $(\tau_k)_{k\geq1}$ and $(\sigma_k)_{k\geq1}$. 
	
	Also, for each $k\geq1$ denote $I_k\equiv\sigma_k-\tau_{k-1}$ and notice that $I_1,I_2,\ldots$ is  a sequence of iid  random variables which have an exponential distribution with rate $\lambda$. Furthermore, for each $k\geq1$, let $N_k$ be the number of jumps of $J(\cdot)$ on $[\sigma_k,\tau_k]$. Note that $N_1,N_2,\ldots$ is a sequence of iid random variables which are distributed according to $N(\cdot)$ (explicitly given in Proposition \ref{prop: PGF}). In a similar fashion, denote the number of jumps of $J(\cdot)$ on $(0,\tau_0]$ by $M$ and notice that $M$ depends on $v$. Furthermore, it is important to notice that the random objects $M$, $(I_k)_{k\geq1}$ and $(N_k)_{k\geq1}$ are all independent. 
	
	 The following identity, which directly follows from the pictorial illustration in Figure 1, is a key ingredient for the rest of our analysis:
	\begin{equation}\label{eq: key identity}
	E_v(x)=xM+\sum_{k=1}^\infty N_k\left(x-\sum_{j=1}^kI_j\right)^+\ \ , \ \ \forall x\geq0\,.
	\end{equation}
	\begin{figure}
	\hspace{-1cm}
	\begin{tikzpicture}
	\draw[thick,->] (0,0) -- (11.5,0) node[anchor=north west]{$\text{Time}$};
	\draw[thick,->] (0,0) -- (0,4) node[anchor=north west]{$\text{Workload}$};
	
	\draw (0,2) coordinate  node [anchor=east] {$v$};
	\draw (0,3) coordinate  node [anchor=east] {$v+x$};

	\draw [line width=0.8pt, blue] (0,2) -- (1,1);
	\draw [{Circle[length=3pt]}-{},line width=0.8pt,blue](1,1.5)--(2.5,0);
	\draw[dashed] (1,0)--(1,2.5);
	\draw (1,0) coordinate  node [anchor=north] {$T_1$};
	\draw [line width=0.8pt, red] (0,3) -- (1,2);
	\draw [{Circle[length=3pt]}-{},line width=0.8pt,red](1,2.5)--(2.5,1);
	
	\draw[line width=0.8pt, blue] (2.5,0) -- (3,0);
	\draw[line width=0.8pt, red] (2.5,1) -- (3,0.5);
	
	\draw[dashed] (3,0)--(3,2.5);
	\draw (3,0) coordinate  node [anchor=north] {$\sigma_1$};
	
	\draw (2.5,0) coordinate  node [anchor=north] {$\tau_0$};

	\draw [{Circle[length=3pt]}-{},line width=0.8pt,blue](3,2)--(4.75,0.25);
	\draw[dashed] (4.75,0)--(4.75,1.5);
	\draw (4.75,0) coordinate  node [anchor=north] {$T_3$};
	\draw [{Circle[length=3pt]}-{},line width=0.8pt,red](3,2.5)--(4.75,0.75);
	\draw [{Circle[length=3pt]}-{},line width=0.8pt,blue](4.75,1)--(5.75,0);
	\draw [{Circle[length=3pt]}-{},line width=0.8pt,red](4.75,1.5)--(6,0.25);
	
	\draw [line width=0.8pt,blue](5.75,0)--(6,0);
	
	\draw[dashed] (6,0)--(6,1);
	\draw (6,0) coordinate  node [anchor=north] {$\ \ \  \sigma_2$};
	\draw (5.75,0) coordinate  node [anchor=north] {$\tau_1$};

	\draw[dashed] (6.25,0)--(6.25,3);
	\draw (6.25,0) coordinate  node [anchor=north] {$\ \ \  \ \ T_5$};
	
	\draw [{Circle[length=3pt]}-{},line width=0.8pt,blue](6,0.75)--(6.25,0.5);
	\draw [{Circle[length=3pt]}-{},line width=0.8pt,red](6,1)--(6.25,0.75);
	
	\draw [{Circle[length=3pt]}-{},line width=0.8pt,blue](6.25,2.75)--(8,1);
	
	\draw [{Circle[length=3pt]}-{},line width=0.8pt,red](6.25,3)--(8,1.25);
	\draw [{Circle[length=3pt]}-{},line width=0.8pt,blue](8,2)--(9.5,0.5);
	\draw [{Circle[length=3pt]}-{},line width=0.8pt,red](8,2.25)--(9.5,0.75);
	\draw [{Circle[length=3pt]}-{},line width=0.8pt,blue](9.5,1.25)--(10.75,0);

	\draw[dashed] (8,0)--(8,2.25);
	\draw (8,0) coordinate  node [anchor=north] {$T_6$};
	
	\draw[dashed] (9.5,0)--(9.5,1.5);
	\draw (9.5,0) coordinate  node [anchor=north] {$T_7$};

	\draw [line width=0.8pt,blue](10.75,0)--(11,0);
	\draw [{Circle[length=3pt]}-{},line width=0.8pt,red](9.5,1.5)--(11,0);
	
	\draw (10.6,0) coordinate  node [anchor=north] {$\ \ \tau_2$};
	\draw (11,0) coordinate  node [anchor=north] {$\ \ \ \zeta$};
	
	\end{tikzpicture}	
	\caption{The blue (resp. red) graph represents the workload process when the initial workload level is $v$ (resp. $v+x$). Note that each jump which occurs during the interval $[0,\tau_0]$ contributes $x$ to the externality. Similarly, each jump which occurs during the interval $[\sigma_1,\tau_1]$ contributes $x-(\sigma_1-\tau_0)$ to the externality. In addition, each jump which occurs during the interval $[\sigma_2,\tau_2]$ contributes $x-(\sigma_1-\tau_0)-(\sigma_2-\tau_1)$ to the externality. Finally, notice that all jumps which occur after the `coupling time' $\zeta$ have no contribution to the externality, and hence the conclusion is that for the current realization we have that $E_v(x)=x+2\left[x-(\sigma_1-\tau_0)\right]+4\left[x-(\sigma_1-\tau_0)-(\sigma_2-\tau_1)\right]$.}
	\end{figure}
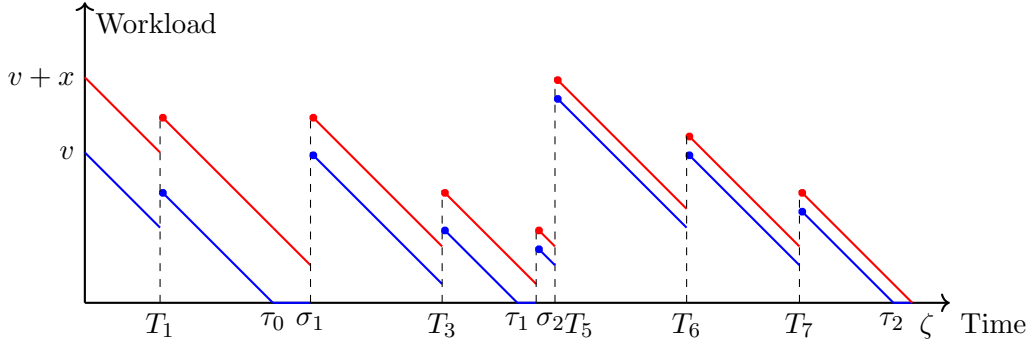
	
	\begin{theorem}\label{thm: marginal}
		For every $x\geq0$ denote, 
		\begin{equation}
		\xi(x)\equiv\min\left\{k\geq1 :\ \sum_{j=1}^kI_j>x\right\}-1\,.
		\end{equation}
		In addition, define a right-continuous nondecreasing stochastic process (in $x$) as follows:
		\begin{equation}
		\dot{E}_v(x)\equiv M+\sum_{k=1}^{\xi(x)}N_k\ \ , \ \ x\geq0\,.
		\end{equation}
		Then, for each $x\geq0$,
		\begin{equation}\label{eq: integral}
		E_v(x)=\int_0^x\dot{E}_v(y)\,\rd y,
		\end{equation}
		and hence $E_v(\cdot)$ is convex with a right-derivative which equals $\dot{E}_v(\cdot)$.
	\end{theorem}
	\begin{remark}
		\normalfont Theorem \ref{thm: marginal} implies that $E_v(\cdot)$ is a {\it convex stochastic process}. For more examples of convex stochastic processes which arise in different applications, see  \cite{Jacobovic2020}.
	\end{remark}
  For each $y\geq0$, let $S(y)$ be the number of jumps that $J(\cdot)$ has until
   \begin{equation}
       \inf\left\{t\geq0:J(t)-t\leq -y\right\}\,.
   \end{equation}
  Notice that $S(v+y)=\dot{E}_v(y)$ for every $y\geq0$. Therefore, when replacing $\dot{E}_v(\cdot)$ by $S(v+\cdot)$ in \eqref{eq: integral}, this equation remains valid. Furthermore, the same technique which was applied in the proof of Proposition \ref{prop: recursive relation} can be used in order to show that $S(\cdot)$ is a compound Poisson process with rate $\lambda$ and jump distribution $N(\cdot)$. As a result, we obtain the following compact representation of the externalities process. 
	
	\begin{corollary}\label{cor: convex}
		In the same probability space in which the model is defined, there is a compound Poisson process $S(\cdot)$ with rate $\lambda$ and jump distribution $N(\cdot)$ such that
		\begin{equation}
		E_v(x)=\int_0^x S(v+y)\,\rd y\ \ , \ \ \forall x\geq0\,.
		\end{equation}
	\end{corollary}
	
	\subsection{Decomposition 2}
	It is interesting to notice that $\dot{E}_v(\cdot)$ equals the number of jumps of $J(\cdot)$ which cause an increase in the value of $E_v(x)$. Consider some arbitrary $x_1,x_2\geq0$ and denote
	\begin{equation}
	\Delta(x_1+x_2,x_1)\equiv E_v(x_1+x_2)-E_v(x_1)\,.
	\end{equation} 
	It is illustrated in Figure 2 that every jump of $J(\cdot)$ which causes an increase in the value of $E_v(x_1)$ contributes $x_2$ to the value of $\Delta(x_1+x_2,x_1)$. This means that we can write
	\begin{equation}
	\Delta_v(x_1+x_2,x_1)=x_2\dot{E}_v(x_1)+\sum_{k=\dot{E}_v(x_1)+1}^\infty N_k\left(x_2-\sum_{j=1}^kI_j\right)^+\,.
	\end{equation}
	Especially, since the workload process is strong Markov, the sum in the right-hand side is distributed as $E_0(x_2)$ and is independent of $\dot{E}_v(x_1)$ (see also Figure 2). 
	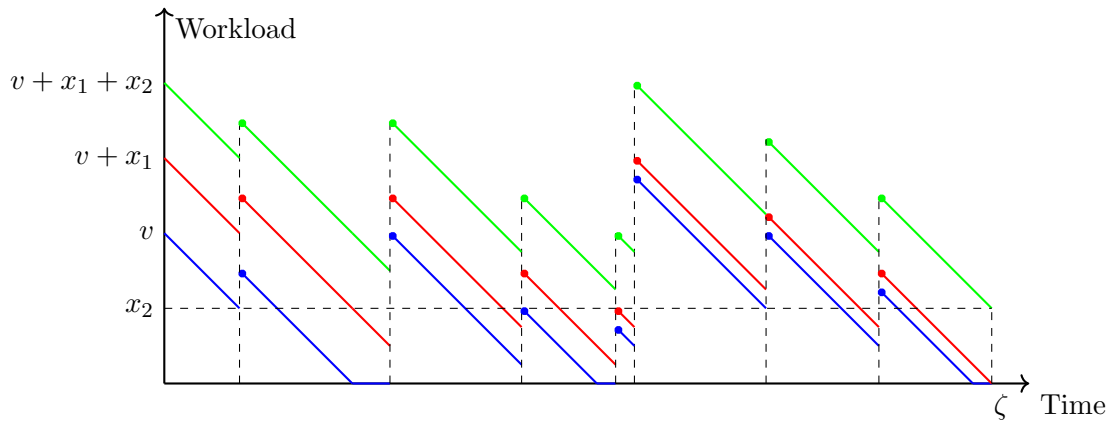
\begin{figure}\label{fig: delta}
		
		\begin{tikzpicture}\hspace{-1.72cm}
		\draw[thick,->] (0,0) -- (11.5,0)node[anchor=north west]{$\text{Time}$};
		\draw[thick,->] (0,0) -- (0,5) node[anchor=north west]{$\text{Workload}$};
		
		\draw (0,2) coordinate  node [anchor=east] {$v$};
		\draw (0,3) coordinate  node [anchor=east] {$v+x_1$};
		\draw (0,4) coordinate  node [anchor=east] {$v+x_1+x_2$};

		\draw [line width=0.8pt, blue] (0,2) -- (1,1);
		\draw [{Circle[length=3pt]}-{},line width=0.8pt,blue](1,1.5)--(2.5,0);
		\draw[dashed] (1,0)--(1,3.5);
		\draw [line width=0.8pt, red] (0,3) -- (1,2);
		\draw [line width=0.8pt, green] (0,4) -- (1,3);
		
		\draw [{Circle[length=3pt]}-{},line width=0.8pt,red](1,2.5)--(2.5,1);
		\draw [{Circle[length=3pt]}-{},line width=0.8pt,green](1,3.5)--(2.5,2);

		\draw[line width=0.8pt, blue] (2.5,0) -- (3,0);
		\draw[line width=0.8pt, red] (2.5,1) -- (3,0.5);
		\draw[line width=0.8pt, green] (2.5,2) -- (3,1.5);

		\draw[dashed] (3,0)--(3,3.5);

		\draw [{Circle[length=3pt]}-{},line width=0.8pt,blue](3,2)--(4.75,0.25);
		\draw[dashed] (4.75,0)--(4.75,2.5);
		\draw [{Circle[length=3pt]}-{},line width=0.8pt,red](3,2.5)--(4.75,0.75);
		\draw [{Circle[length=3pt]}-{},line width=0.8pt,green](3,3.5)--(4.75,1.75);
		
		\draw [{Circle[length=3pt]}-{},line width=0.8pt,blue](4.75,1)--(5.75,0);
		\draw [{Circle[length=3pt]}-{},line width=0.8pt,red](4.75,1.5)--(6,0.25);
		\draw [{Circle[length=3pt]}-{},line width=0.8pt,green](4.75,2.5)--(6,1.25);

		\draw [line width=0.8pt,blue](5.75,0)--(6,0);
		
		\draw[dashed] (6,0)--(6,2);

		\draw[dashed] (6.25,0)--(6.25,4);
		
		\draw [{Circle[length=3pt]}-{},line width=0.8pt,blue](6,0.75)--(6.25,0.5);
		\draw [{Circle[length=3pt]}-{},line width=0.8pt,red](6,1)--(6.25,0.75);
		\draw [{Circle[length=3pt]}-{},line width=0.8pt,green](6,2)--(6.25,1.75);

		\draw [{Circle[length=3pt]}-{},line width=0.8pt,blue](6.25,2.75)--(8,1);
		
		\draw [{Circle[length=3pt]}-{},line width=0.8pt,red](6.25,3)--(8,1.25);
		\draw [{Circle[length=3pt]}-{},line width=0.8pt,green](6.25,4)--(8,2.25);
		
		\draw [{Circle[length=3pt]}-{},line width=0.8pt,blue](8,2)--(9.5,0.5);
		\draw [{Circle[length=3pt]}-{},line width=0.8pt,red](8,2.25)--(9.5,0.75);
		\draw [{Circle[length=3pt]}-{},line width=0.8pt,green](8,3.25)--(9.5,1.75);
		
		\draw [{Circle[length=3pt]}-{},line width=0.8pt,blue](9.5,1.25)--(10.75,0);

		\draw[dashed] (8,0)--(8,3.25);

		\draw[dashed] (9.5,0)--(9.5,2.5);

		\draw [line width=0.8pt,blue](10.75,0)--(11,0);
		\draw [{Circle[length=3pt]}-{},line width=0.8pt,red](9.5,1.5)--(11,0);
		\draw [{Circle[length=3pt]}-{},line width=0.8pt,green](9.5,2.5)--(11,1);

		\draw (11,0) coordinate  node [anchor=north] {$\ \ \zeta$};
		\draw[dashed] (11,0)--(11,1);
		\draw[dashed] (0,1)--(11,1);
		\draw (0,1) coordinate  node [anchor=east] {$x_2$};
		\end{tikzpicture}\caption{The green (resp.\ red) graph describes a sample path of the workload process when a customer $c$ with a service demand of $x_1+x_2>0$ (resp. $x_1$) arrives at time zero and sees a system with existing workload level $v>0$. The blue graph describes a sample path of the workload of the same system once $c$ reduces her service requirement to zero. Note that the jumps of the graphs are coordinated. In fact, each jump is associated with an arrival of a customer and the size of the jump is the service demand of that customer. Observe that every jump on $(0,\zeta)$ adds $x_2$ to $\Delta_v(x_1+x_2,x_1)$. In addition, the value of the green graph at $\zeta$ equals $x_2$. Thus, a regenerative argument yields that $\Delta_v(x_1+x_2,x_1)$ is distributed as $x_2$ multiplied by the number of jumps on $(0,\zeta)$ plus an independent random variable which is distributed like $E_0(x_2)$.}
	\end{figure}
	
	Furthermore, assume that $\xi\sim\text{Poi}(\lambda x_2)$ and $U_1,U_2,\ldots$ is an iid sequence of random variables which are distributed uniformly on $[0,1]$. In particular, assume that $\xi$, $(U_k)_{k\geq1}$ and $(N_k)_{k\geq1}$ are independent. In addition, for each $j\geq1$ we use the notation
	\begin{equation}
	U_{(1),j}\geq U_{(2),j}\geq\ldots\geq U_{(j),j}
	\end{equation}
	in order to denote the order statistics of $U_1,U_2,\ldots,U_j$. Then an application of known `symmetry properties' yields the following distributional equality:
	\begin{equation}
	E_0(x_2)=\sum_{k=1}^\infty N_k\left(x_2-\sum_{j=1}^kI_j\right)^+\stackrel{d}{=}\sum_{k=1}^{\xi}N_k\left(x_2U_{(1),\xi}\right)\stackrel{d}{=}x_2\sum_{k=1}^{\xi}N_kU_k\,.
	\end{equation}
	This argument can be applied recursively in order to derive the following theorem.   As illustrated in Section \ref{sec: moments}, it provides us with a systematic approach to compute the moments of the finite-dimensional distributions of $E_v(\cdot)$. 
	\begin{theorem}\label{thm: decomposition}
		Let $x_1,x_2,\ldots,x_k\geq0$ and assume that $\mathcal{N}\equiv \left\{N_{m,i}:i,m\geq1\right\}$, $\mathcal{U}\equiv\left\{U_{m,i}:i,m\geq1\right\}$ and $\left(\xi_j\right)_{j\geq1}$ are such that:
		\begin{enumerate}
			\item $\mathcal{N}$ is an infinite array of iid random variables such that $N_{1,1}$ is distributed according to $N(\cdot)$.
			
			\item $\mathcal{U}$ is an infinite array of iid random variables which are distributed uniformly on $[0,1]$. 
			
			\item $\xi_1,\xi_2,\ldots$ are independent random variables such that $\xi_1\sim\text{\rm Poi}(\lambda v)$ and $\xi_j\sim\text{\rm Poi}(\lambda x_{j-1})$, $j\geq2$.  
			
			\item $\mathcal{N},\mathcal{U}$ and $(\xi_j)_{j\geq1}$ are independent.  
		\end{enumerate}
	Then, 
	\begin{equation}\label{eq: decomposition 1}
	E_v\left(\sum_{i=1}^kx_i\right)\stackrel{d}{=}\sum_{j=1}^kx_j\left(\sum_{l=1}^j\sum_{m=1}^{\xi_l}N_{m,l}+\sum_{m=1}^{\xi_{j+1}}N_{m,j+1}U_{m,j+1}\right)
	\end{equation}
	and 
	\begin{equation}\label{eq: decomposition 2}
	\Delta_v(x_1+x_2,x_1)\stackrel{d}{=}x_2\left(\sum_{m=1}^{\xi_1}N_{m,1}+\sum_{m=1}^{\xi_2}N_{m,2}+\sum_{m=1}^{\xi_{3}}N_{m,3}U_{m,3}\right)\,.
	\end{equation}
	\end{theorem}
	
	\section{Moments of the finite-dimensional distributions}\label{sec: moments}
	This section concentrates on the evaluation of moments corresponding to the finite dimensional distributions of the externalities process $E_v(\cdot)$. We first present the mean and variance, then the auto-covariance and auto-correlation, after which we proceed with higher moments.
	\subsection{Mean and variance}
	Fix $x>0$ and notice that an insertion of $k=1$ into \eqref{eq: decomposition 1} yields that
	\begin{equation}\label{eq: univariate decomposition}
		E_v(x)\stackrel{d}{=}x\left(\sum_{m=1}^{M_1}N_{m,1}+\sum_{m=1}^{M_{2}}N_{m,2}U_{m,2}\right)\,.
	\end{equation}
	Thus, by an application of the formula of an expectation of a compound Poisson random variable, we directly obtain that
	\begin{equation}\label{eq: expectation}
		\mathbb{E}[E_v(x)]=\lambda x\left(v+\frac{x}{2}\right) \eta_1
	\end{equation}
	and, as expected, $\mathbb{E}[E_v(x)]$ is convex in $x$. Similarly, the formula of the variance of a compound Poisson random variable may be used in order to derive that
	\begin{equation}\label{eq: variance}
		\text{Var}\left[E_v(x)\right]=\lambda x^2\left(v+\frac{x}{3}\right)\eta_2\,.
	\end{equation}
	
	\subsection{Auto-covariance and auto-correlation}
	Fix some $x_1,x_2>0$.
	Since the sums in the right-hand side of \eqref{eq: decomposition 2} are independent, we find 
	\begin{equation}
		\text{Var}\left[\Delta_v(x_1,x_1+x_2)\right]=\lambda x_2^2\left(v+x_1+\frac{x_2}{3}\right)\eta_2\,.
	\end{equation} 
	In addition,
	\begin{align}
	\text{Var}\left[\Delta_v(x_1,x_1+x_2)\right]&=\text{Var}\left[E_v(x_1+x_2)\right]\\&+\text{Var}\left[E_v(x_1)\right]-2\text{Cov}\left[E_v(x_1+x_2),E_v(x_1)\right]\,,\nonumber
	\end{align}
	and hence an insertion of \eqref{eq: variance} implies that the auto-covariance function of $E_v(\cdot)$ equals 
	\begin{align}\label{eq: covariance}
	R_v(x_1,x_1+x_2)&\equiv\text{Cov}\left[E_v(x_1+x_2),E_v(x_1)\right]\\&=\frac{\lambda\eta_2}{2}\cdot\frac{2x_1^3+6vx_1^2+3x_1^2x_2+6vx_1x_2}{3}\,.\nonumber
	\end{align}
As argued in the introduction, in the situation of a customer arriving at time 0 with {\it two} tasks (of size $x_1$ and $x_2$, respectively), $E_v(x_1+x_2)-E_v(x_1)$ represents the total waiting time that could be saved by the other customers if the customer gave up on their second task but insisted on completing the first one. The auto-covariance \eqref{eq: covariance} provides insight into the effect of the additional $x_2$.

	As a result, the auto-correlation function is given by
	\begin{align}\label{eq: correlation}
		\rho_v(x_1,x_1+x_2)&\equiv\text{Corr}\left[E_v(x_1+x_2),E_v(x_1)\right]\\&=\frac{2x_1^3+6vx_1^2+3x_1^2x_2+6vx_1x_2}{6x_1(x_1+x_2)\sqrt{\left(v+\frac{x_1}{3}\right)\left(v+\frac{x_1+x_2}{3}\right)}}\,.\nonumber
	\end{align}
	Surprisingly, the expression in \eqref{eq: correlation} is {\it invariant} with respect to the service distribution and the arrival rate. At the same time,  observe that $R_v(x_1,x_1+x_2)$ is positive. In addition, the expression of $R_v(x_1,x_1+x_2)$ actually shows that the externalities process is not {\it wide sense stationary} (see the definition in  \cite[p.~15]{Yaglom2004}).
	
	\begin{remark}
	    \normalfont Later, in Section \ref{sec: v is random} we consider a setup in which $v$ is a general nonnegative random variable, independent from the arrival process and service requirements. There, it is shown that in the more complex setup, the auto-correlation function depends on the arrival rate and service distribution unless $v$ is a degenerate random variable. 
	\end{remark}

	\subsection{Higher moments}
	Higher moments (including joint moments) of $E_v(\cdot)$ may be derived via differentiation of the LST formula which is given in the next theorem. This is a tedious derivation that we decided to leave out.   The below result is particularly useful when analyzing a situation in which the customer arriving at time 0 has $k$ tasks, having sizes $x_1,\ldots,x_k$.
	
	\begin{theorem}\label{thm: LST}
		Let $k\geq1$ and $\alpha\equiv(\alpha_1,\alpha_2,\ldots,\alpha_k),x\equiv(x_1,x_2,\ldots,x_k)\in(0,\infty)^k$. In addition, define
		\begin{equation}
			w(x,\alpha)\equiv\alpha_1x_1+\alpha_2(x_1+x_2)+\ldots+\alpha_k\sum_{i=1}^kx_i\,,
		\end{equation}
		and for every $u\in(0,1)$ denote
		\begin{align}
			&s_1(x,\alpha,u)\equiv\alpha_1x_1u+\alpha_2(x_1u+x_2)+\ldots+\alpha_k\left(x_1u+\sum_{i=2}^kx_i\right)\,,\\&s_2(x,\alpha,u)\equiv\alpha_2x_2u+\alpha_3(x_2u+x_3)+\ldots+\alpha_k\left(x_2u+\sum_{i=3}^kx_i\right)\,,\nonumber\\&\ \ \ \ \ \ \ \ \ \ \ \ \ \ \vdots\nonumber\\&s_k(x,\alpha,u)\equiv\alpha_kx_ku\,.
		\end{align}
		Then, for any $v>0$,
		\begin{align}
			\mathbb{E}\exp\left\{-\sum_{l=1}^k\alpha_lE_v\left(\sum_{i=1}^lx_i\right)\right\}&=\exp\left\{\lambda v\left[\widetilde{N}\left(w(x,\alpha)\right)-1\right]\right\}\\&\cdot\prod_{l=1}^{k}\int_0^1\exp\left\{\lambda x_l\left[\widetilde{N}\left(s_l(x,\alpha,u)\right)-1\right]\right\}\rd u\,.\nonumber
		\end{align}
	\end{theorem}
	
	\section{Crossing times of $\dot{E}_0(\cdot)$}\label{sec: crossing times}
		The process $\dot{E}_v(\cdot)$ is nondecreasing such that $\dot{E}_v(0)=M$, and $E_v(x)\uparrow\infty$ as $x\to\infty$. Therefore, it is natural to study the crossing times of the process $\dot{E}_v(\cdot)$. Namely, fix some $y>0$ and the corresponding crossing time is  
		\begin{equation}
			x_v(y)\equiv\inf\left\{x\geq0:\dot{E}_v(x)\geq y\right\}\,.
		\end{equation}
		In the sequel we consider the special case $v=0$ for which $M=0$ and hence $x(y)\equiv x_0(y)$ has a relatively tractable representation (but see Remark \ref{vpos} below for some reflections on the case $v>0$). In fact, Theorem \ref{thm: decomposition} yields that 
		\begin{align}\label{eq: crossing}
		x(y)&\equiv\inf\left\{x\geq0:\dot{E}_0(x)\geq y\right\}\\&=\inf\left\{x\geq0:\sum_{m=1}^{\xi(x)}N_k\geq \lceil y\rceil\right\}=\sum_{k=1}^{\upsilon(y)}I_k\,,\nonumber
		\end{align}
		where $\upsilon(y)\equiv\min\left\{t\geq1:\sum_{m=1}^{t}N_k\geq\lceil y\rceil\right\}$. In particular, notice that $\upsilon(y)$ and $I_1,I_2,\ldots$ are independent. Furthermore, $\upsilon(y)$ can be described as the time until absorption in a Markov chain with a unique absorbing state. Specifically, this chain has a state-space $\left\{0,1,2,\ldots,\lceil y\rceil\right\}$ with an absorbing state $\lceil y\rceil$ and an initial state $0$. In addition, the transition matrix is equal to
		$P\equiv\left[p_{ij}\right]_{1\leq i,j\leq\lceil y\rceil}$, given by
		\begin{equation}
			P=\begin{pmatrix}
				0 & N(1) & N(2)& N(3)&\cdots & N\left(\lceil y\rceil-1\right) & 1-\sum_{s=1}^{\lceil y\rceil-1}N(s) \\
				0 & 0 & N(1) & N(2)& \cdots & N\left(\lceil y\rceil-2\right)& 1-\sum_{s=1}^{\lceil y\rceil-2}N(s) \\
				0 & 0 & 0& N(1) & \ldots & N\left(\lceil y\rceil-3\right) & 1-\sum_{s=1}^{\lceil y\rceil-3}N(s)\\
				\vdots  & \vdots  & \vdots& \vdots& \ddots & \vdots & \vdots  \\
				0 & 0 & 0 & 0 &\cdots & N(1) & 1-N(1)\\0 & 0 &0 &0& \ldots& 0 &1
			\end{pmatrix}\,.
		\end{equation}
		It is well-known, that the mean of $\upsilon(y)$ can be characterized via
		\begin{equation}
			\psi_0(y)\equiv\mathbb{E}\upsilon(y)=1+\sum_{k=1}^{\lceil y\rceil-1}p_{0k}\psi_k=1+\sum_{k=1}^{\lceil y\rceil-1}N(k)\psi_k \,,
		\end{equation}
		where $\psi_{\lceil y\rceil-1}=1$ and $\psi_1,\psi_2,\ldots,\psi_{\lceil y\rceil-2}$ are given recursively by the equations
		\begin{equation}
			\psi_k=1+\sum_{i=k+1}^{\lceil y\rceil-1}p_{ki}\psi_i=1+\sum_{i=k+1}^{\lceil y\rceil-1}N(i-k)\psi_i\,.\ \ , \ \ k=1,2,\ldots,\lceil y\rceil-2\,.
		\end{equation}
		Therefore, Wald's identity may be applied to \eqref{eq: crossing} to deduce that
		\begin{equation}\label{eq: mean crossing}
			\mathbb{E}x(y)=\mathbb{E}\upsilon(y)\,\mathbb{E}I_{2,1}=\frac{\psi_0(y)}{\lambda}\,.
		\end{equation}
		
		\begin{remark}
			\normalfont The second moment of $\upsilon(y)$ can be computed using a similar technique, thus also yielding $\text{Var}\left[\upsilon(y)\right]$. Hence, it is possible to compute the variance of $x(y)$ via the formula
			\begin{equation}
				\text{Var}\left[x(y)\right]=\frac{\psi_0(y)+\text{Var}\left[\upsilon(y)\right]}{\lambda^2}\,.
			\end{equation}
		\end{remark}
		\begin{remark} \normalfont \label{vpos}
			When $v>0$, it makes sense to rely on a similar computation in which we condition and de-condition on $M$. In practice, we do not see how this computation leads to a tractable expression for the general case. 
		\end{remark}

	\section{Gaussian approximation of $E_v(\cdot)$}\label{sec: heavy traffic}
	The main result of this section concerns a Gaussian approximation for the externalities process. In order to provide an accurate statement of this result, notice that the model which was described in Section \ref{sec: definition} is characterized by the triplet $(\lambda,B(\cdot),v)$. Fix $v\geq0$ and consider a {\it sequence} of models
	\begin{equation}
	(\lambda_n,B_n(\cdot),v)\ \ , \ \ n\geq1\,,
	\end{equation} such that the $n$-th model is associated with an arrival rate $\lambda_n>0$ and a service distribution $B_n(\cdot)$. Respectively, for each $n\geq1$, we introduce the notation
	\begin{align}
	&\mu_{k,n}\equiv\int_0^\infty t^k\rd B_n(t)\ \ , \ \ k\geq1\,,\\&
	\rho_n\equiv\lambda_n\mu_{1,n}\,.\nonumber
	\end{align}
	In addition, for each $n\geq1$, denote the externalities process which is associated with the $n$-th model by $E_v^{(n)}(\cdot)$. 
	Also, let $N_n(\cdot)$ be the probability mass function of the number of customers who got service during a single busy period of a FCFS M/G/1  queue with an arrival rate $\lambda_n$ and a service distribution $B_n(\cdot)$. Correspondingly, for each $n,k\geq1$ denote
	\begin{equation}
		\eta_{k,n}\equiv\sum_{s=1}^\infty s^kN_n(s)\ \ , \ \ k\geq1,
	\end{equation} 
	and observe that $\eta_{k,n}$ is the analogue of $\eta_k$ in the $n$-th model. 
	
	\subsection{Functional central limit theorem}
	The main result of this section is stated in the next functional central limit theorem. 
	\begin{theorem}\label{thm: externalities - Gaussian approx} 
		Define, for a fixed $v\ge0$,
		\begin{enumerate}
			\item A stochastic process
			\begin{equation}
				H_v(x)\equiv \int_0^xW(v+y)\,\rd y\ \ , \ \ x\geq0
			\end{equation}
			such that $W(\cdot)$ is a standard Wiener process.
			
			\item A sequence (in $n=1,2,\ldots$) of stochastic processes 
			\begin{equation}
				\hat{E}_v^{(n)}(x)\equiv
				\frac{E^{(n)}_v(x)-\mathbb{E}\big[E^{(n)}_v(x)\big]}{\sqrt{\lambda_n\eta_{2,n}}} \ \ , \ \ x\geq0\,.
			\end{equation}
		\end{enumerate} 
		In addition, assume that the next conditions hold:
		\begin{description}
			\item [(i)] $\lambda_n\rightarrow\infty$ as $n\to\infty$.
			
			\item[(ii)] There is $n'\geq1$ such that $\rho_n<1$ for every $n\geq n'$.
			
			\item [(iii)] $\frac{\eta_{3,n}}{\sqrt{\lambda_n\eta^3_{2,n}}}\rightarrow0$ as $n\to\infty$.
		\end{description}  
	Then,   
	\begin{equation}\label{eq: convergence}
	\hat{E}_v^{(n)}(\cdot)\Rightarrow H_v(\cdot)\ \ \text{as}\ \ n\to\infty
	\end{equation}
	where $\Rightarrow$ denotes weak convergence on $\mathcal{D}[0,\infty)$ equipped with the uniform metric {(on compacta)}.
	\end{theorem}
	Observe that checking Condition \textbf{(iii)} is not straightforward because it is phrased in terms of the moments of $N_n(\cdot)$. The following proposition presents two sets of sufficient conditions which are considerably more easy to verify. Broadly speaking, the proof of these sets of conditions being sufficient relies on the expressions appearing in the statement of Corollary~\ref{cor: N moments}.
	\begin{proposition}\label{prop: sufficient conditions}
	
	   $\text{ }$\newline{} \begin{enumerate}
	        \item Assume that Condition \textbf{(i)} and in addition the two conditions    \begin{equation}\label{eq: sufficient conditions 1}
	            \limsup_{n\to\infty}\rho_n<1\ , \  \limsup_{n\to\infty}\lambda_n^3\mu_{3,n}<\infty\,,
	        \end{equation}
	       are all satisfied. Then, \eqref{eq: convergence} is valid. 
	        
	        \item Assume that Condition \textbf{(i)}, Condition \textbf{(ii)} and in addition the  three conditions    \begin{equation}\label{eq: sufficient conditions 2}
	            \liminf_{n\to\infty}\lambda_n^2\mu_{2,n}>0\ , \  \limsup_{n\to\infty}\lambda_n^3\mu_{3,n}<\infty\ ,\ \lim_{n\to\infty}\lambda_n(1-\rho_n)=\infty\,,
	        \end{equation}
	       are all satisfied. Then, \eqref{eq: convergence} is valid.  
	    \end{enumerate}
	\end{proposition}
	
	\begin{remark}
	    \normalfont Note that the condition
	    \begin{equation*}
	        \limsup_{n\to\infty}\rho_n<1
	    \end{equation*}
	   implies Condition \textbf{(ii)}. In addition, it does not go together with a heavy-traffic regime (i.e., $\rho_n\uparrow1$ as $n\to\infty$) under the first set of conditions in Proposition \ref{prop: sufficient conditions}.
	   Thus, the added value of the second set of conditions in Proposition \ref{prop: sufficient conditions} is that it  could cover a heavy-traffic regime, i.e., $\rho_n\uparrow1$ as $n\to\infty$.
	\end{remark}
	
	\begin{remark}
	    \normalfont Assume that Condition \textbf{(i)} is satisfied. Thus, if
	    \begin{equation}\label{eq: condi1}
	    \lim_{n\to\infty}\sqrt{\lambda_n}\left(1-\rho_n\right)=c    
	    \end{equation}
	    for some $c\in(0,\infty)$, then
	    \begin{equation}\label{eq: condi2}
	    \lim_{n\to\infty}\lambda_n\left(1-\rho_n\right)=\infty\,.
	    \end{equation}
	    At the same time, observe that \eqref{eq: condi1} is not necessary for \eqref{eq: condi2} even under the assumption that Condition \textbf{(i)} is satisfied. 
	\end{remark}
	
	The general idea of the proof of Theorem \ref{thm: externalities - Gaussian approx} is as follows: Notice that Condition \textbf{(ii)} allows us to apply Corollary \ref{cor: convex} and hence  for each $n\geq n'$ there is a \textit{compensated} compound Poisson process $S_n(\cdot)$ with rate $(\lambda_n)$ and jump distribution $(N_n(\cdot))$ such that
	\begin{equation}
	\hat{E}_n^{(n)}(x)=\int_0^x\frac{S_n(v+y)}{\sqrt{\lambda_n\eta_{2,n}}}\rd y\ \ , \ \ \forall x\geq0\,.
	\end{equation}   
	Then, the crucial part of the proof is to show that
	\begin{equation}\label{eq: benchmark}
	\frac{S_n(\cdot)}{\sqrt{\lambda_n\eta_{2,n}}}\Rightarrow W(\cdot)\ \ \text{as}\ \ n\to\infty\,;
	\end{equation}   
	the rest will follow from this benchmark via standard arguments. In the upcoming Section \ref{subsec: CPP Gaussian approx} we address a general result about a Gaussian approximation of compound Poisson processes. This will help in proving \eqref{eq: benchmark}. \begin{remark}
	    \normalfont An extensive account of heavy-traffic approximations of queueing systems can be found in \cite{Whitt2002}. Notably, heavy-traffic approximations have been developed for various functionals of the queueing process (such as the number of customers and the waiting time), but to the best of our knowledge we are the first to do so for the externalities process. This means that, in the strict sense, we cannot compare our Theorem \ref{thm: externalities - Gaussian approx}  with existing results. 
This being said, there is a vast literature on Gaussian approximations for sequences of compound Poisson processes, related to Theorem~\ref{thm: CPP-Gaussian approx} below (which is heavily relied upon in our derivation of Theorem~\ref{thm: externalities - Gaussian approx}); we therefore include in Section \ref{subsec: CPP Gaussian approx}  a comparison between Theorem~\ref{thm: CPP-Gaussian approx} and related results. 
	\end{remark}
	
	We proceed by discussing an immediate implication of Theorem \ref{thm: externalities - Gaussian approx}. To this end, fixing $x\geq0$, recall that it is well-known result that 
	\begin{equation}\label{eq: Ito}
	\int_0^xW(v+y)\,\rd y\sim\mathcal{N}\left(0,x^2v+\frac{x^3}{3}\right)\,.
	\end{equation}
	Hence, under the conditions of Theorem \ref{thm: externalities - Gaussian approx} we conclude the following convergence:
	\begin{equation}
	\hat{E}_v^{(n)}(x)\xrightarrow{d}\mathcal{N}\left(0,x^2v+\frac{x^3}{3}\right)\ \ \text{as}\ \ n\to\infty\,.
	\end{equation}
	
	\begin{remark}
		\normalfont In fact, taking into account \eqref{eq: variance}, the current analysis gives a new proof for \eqref{eq: Ito} which is not based on stochastic calculus at all but only on approximation of a standard Wiener process by normalized compensated compound Poisson processes. Since \eqref{eq: Ito} is known and the current proof is not simpler than the existing one, we mention this result in passing. 
	\end{remark}
	
	\subsection{Gaussian approximation to compound Poisson process}\label{subsec: CPP Gaussian approx}
 In this subsection we discuss a general Gaussian approximation result for compound Poisson processes and relate it to existing results. As mentioned, it is used in the proof of Theorem \ref{thm: externalities - Gaussian approx}, but may have broader applications. 

 \subsubsection{Gaussian approximation result}
	The following theorem, proven in Section \ref{proofTh5}, includes a statement about a Gaussian approximation of a general compound Poisson process. Possibly, this theorem may have other applications besides those that appear in the current work. 
	\begin{theorem}\label{thm: CPP-Gaussian approx}
		For each $n\geq1$ let $\left\{J_n(t):t\geq0\right\}$ be a compensated compound Poisson process with rate $\lambda_n>0$ and jump distribution $F_n(\cdot)$ such that
		\begin{equation}
		\sigma_n\equiv\int_{-\infty}^\infty t^2\rd F_n(t)\ \in(0,\infty)\,.
		\end{equation}
		In addition, denote 
		\begin{equation}
			\  \ \ \nu_n\equiv\int_{-\infty}^\infty |t|^3\rd F_n(t)\,,
		\end{equation}
		 and assume that both of the following conditions hold:
		 \begin{description}
		 	\item[(I)] $\lambda_n\rightarrow\infty$ as $n\to\infty$.
		 	
		 	\item [(II)] $\frac{\nu_n}{\sqrt{\lambda_n\sigma_n^3}}\rightarrow0$ as $n\to\infty$. 
		 \end{description}  
	 Then, 
	\begin{equation}
		\frac{J_n(\cdot)}{\sqrt{\lambda_n\sigma_n}}\Rightarrow W(\cdot)\ \ \text{as}\ \ n\to\infty\,,
	\end{equation}
	where $W(\cdot)$ is a standard Wiener process and $\Rightarrow$ denotes weak convergence on $\mathcal{D}[0,\infty)$ equipped with the uniform metric {(on compacta)}.
	\end{theorem}
	\begin{remark}
	    \normalfont Intuitively speaking, Condition \textbf{(I)} implies that the jumps become more frequent as $n\to\infty$ while Condition \textbf{(II)}  makes sure that the jump distribution should not become too `wild' as $n\to\infty$. 
	\end{remark}
	\subsubsection{Comparison with the existing literature}
	Let $J(\cdot)$ be a compensated-compound Poisson process with rate $\lambda$ and a jump distribution with finite fourth moment. Denote the standard deviation of the jump distribution by $\gamma$. Then, \cite[Corollary 3.7]{Khoshnevisan1993} states conditions under which the sequence (in $n$) of processes 
	\begin{equation}
	   \widetilde{J}_n(t)\equiv \frac{J(nt)}{\gamma\sqrt{n}}\ \ , \ \ t\geq0\,,  
	\end{equation}
    weakly converges to a standard Wiener process in $\mathcal{D}[0,1]$ equipped with the Skorohod topology. Thus, some differences between Theorem \ref{thm: CPP-Gaussian approx} and \cite[Corollary 3.7]{Khoshnevisan1993} are:
    \begin{enumerate}
        \item Consider the setup of Theorem \ref{thm: CPP-Gaussian approx} with (i.e.,  $\lambda_n=n\lambda$) and $\sigma_n\equiv \frac{\gamma^2}{\lambda}$. Then, we get that $J_n\equiv\widetilde{J}_n$ and hence in that sense the setup of Theorem \ref{thm: CPP-Gaussian approx} is more general. 
        
        \item Theorem \ref{thm: CPP-Gaussian approx} guarantees weak convergence in a  different topological space.
        
        \item \cite[Corollary 3.7]{Khoshnevisan1993} requires that the fourth moment of the jump distribution is finite while Theorem \ref{thm: CPP-Gaussian approx} imposes no conditions on the fourth moment of $F_n$ (for any $n\geq1$).  
        
        \item In \cite[Corollary 3.7]{Khoshnevisan1993}, we get that $\lambda_n$ grows linearly in $n$ which implies Condition \textbf{(I)}, but obviously Condition \textbf{(I)} might be satisfied in other asymptotic regimes.
        
        \item In \cite[Corollary 3.7]{Khoshnevisan1993}, we get that $\nu_n$ and $\sigma_n$ remain fixed (in $n$) and hence, due to the linear growth of $\lambda_n$, Condition \textbf{(II)} is satisfied. Once again, obviously it could be satisfied in other asymptotic regimes as well.  
    \end{enumerate}
    Another result \cite[Theorem 1.1]{Pang2017}, is about a weak convergence in $\mathcal{D}[0,\infty)$ equipped with the Skorohod topology of a sequence of modulated compound Poisson processes. When all the processes in this sequence are compound Poisson processes (i.e., when the modulating Markov chains in the background are all degenerate ones), then \cite[Assumption 1]{Pang2017} is reduced to:
    \begin{enumerate}
        \item Linear growth of the sequence $\lambda_n$ as $n\to\infty$.  
        
        \item The sequences (in $n$) of the means and standard deviations of $F_n$ should both  converge to constants.  
    \end{enumerate}   
    We conclude that there is a strong resemblance between the comparison of Theorem \ref{thm: CPP-Gaussian approx} with \cite[Theorem 1.1]{Pang2017} and the comparison of Theorem \ref{thm: CPP-Gaussian approx} with \cite[Corollary 3.7]{Khoshnevisan1993}. 
    
    Another strand of literature regards the properties of a sequence of compound Poisson processes which weakly converges to a limiting process (see, e.g., \cite{Lambert2013,Lambert2015,Sarkar2005}). This literature predominantly focuses on necessary conditions for weak convergence of such sequences, while Theorem \ref{thm: CPP-Gaussian approx} presents sufficient conditions. 
    
	\section{When $v$ is a random variable}\label{sec: v is random}
	In this part we revisit some results of the previous sections in the situation that $v$ is a nonnegative random variable which is independent from the arrival process and the service requirements of the customers. The motivation for this extension of the existing framework lies in the fact that if $v$ has the stationary distribution of an M/G/1  queue with an arrival rate $\lambda$ and a service distribution $B(\cdot)$, then we recover the setup of Haviv and Ritov \cite{Haviv1998}. For simplicity of notation, denote the conditional expectation (resp.\ covariance) given $v$ by $\mathbb{E}_v$ (resp. $\text{Cov}_v$).
	
	\subsection{Expressions for moments}
	To begin with, it is immediate that the decompositions of  Section~\ref{sec: definition} remain true  when the initial workload $v$ is a general random variable. Similarly, Theorem \ref{thm: externalities - Gaussian approx} may be phrased in the extended setup. This is because for every bounded uniformly-continuous functional $f$ we may apply the law of total expectation and then apply the dominated convergence theorem with Theorem \ref{thm: externalities - Gaussian approx} so as to deduce the needed result (see also \cite[Corollary IV.9]{Pollard2012}). 
	
	A similar approach may be applied in order to derive the moments of the externalities process. For example, for every $x\geq0$,
	\begin{equation}\label{eq: expectation1}
	\mathbb{E}[E_v(x)]=\mathbb{E}\big[\mathbb{E}_v[E_v(x)]\big]=\lambda x\left(\mathbb{E}v+\frac{x}{2}\right)\eta_1\,.
	\end{equation}  
	Furthermore, for every $x_1,x_2\geq0$ deduce that
	\begin{align}
	\mathbb{E}[\text{Cov}_v[E_v(x_1),\,&E_v(x_1+x_2)]]\notag\\
	&=\mathbb{E}[R_v(x_1,x_1+x_2)]\\&=\frac{\lambda\eta_2}{2}\cdot\frac{2x_1^3+6x_1^2\mathbb{E}v+3x_1^2x_2+6x_1x_2\mathbb{E}v}{3}\nonumber
	\end{align}
	and once $\mathbb{E}v^2<\infty$ we also get
	\begin{align}
	\text{Cov}[\mathbb{E}_v[E_v(x_1)]&,\mathbb{E}_v[E_v(x_1+x_2)]]\notag\\
	&=\text{Cov}\left[\lambda x_1\left(v+\frac{x_1}{2}\right) \eta_1,\lambda (x_1+x_2)\left(v+\frac{x_1+x_2}{2}\right) \eta_1\right]\nonumber\\&=\lambda^2\eta_1^2x_1(x_1+x_2)\text{Var}(v)\,.
	\end{align}
	Thus, the law of total covariance yields that
	\begin{align}\label{eq: covariance1}
	\text{Cov}&\left[E_v(x_1),E_v(x_1+x_2)\right]\\&=\frac{\lambda\eta_2}{2}\cdot\frac{2x_1^3+6x_1^2\mathbb{E}v+3x_1^2x_2+6x_1x_2\mathbb{E}v}{3}+\lambda^2\eta_1^2x_1(x_1+x_2)\text{Var}(v)\,.\nonumber
	\end{align}
	In particular, when $x_1=x$ and $x_2=0$ the variance is obtained:
	\begin{equation}\label{eq: variance1}
	\text{Var}\left[E_v(x)\right]=\left(\lambda x \eta_1\right)^2\text{Var}(v)+\lambda\eta_2 x^2\left(\mathbb{E}v+\frac{x}{3}\right)\,.
	\end{equation}
	\begin{remark}
		\normalfont Equations \eqref{eq: covariance1} and \eqref{eq: variance1} imply that  the correlation is invariant to the arrival rate and the service distribution if and only if $\text{Var}(v)=0$  (or equivalently, when $v$ equals a constant with probability one).  
	\end{remark}
	For higher moments, it is possible to differentiate the LST formula, as  given in the next corollary. Just like in Section \ref{sec: moments}, we do not include these computations here. The proof follows from conditioning and de-conditioning on $v$ with the result of Theorem \ref{thm: LST}.  
	
	\begin{corollary}\label{cor: LST}		
		Let $k\geq1$ and $\alpha\equiv(\alpha_1,\alpha_2,\ldots,\alpha_k),x\equiv(x_1,x_2,\ldots,x_k)\in(0,\infty)^k$. In addition, denote the LST of $v$ by 
		\begin{equation}
		\widetilde{v}(t)\equiv {\mathbb E}\,e^{-tv}\ \ , \ \ t>0
		\end{equation}
		and consider the notation of Theorem \ref{thm: LST}.
		Then, 
		\begin{align}
		\mathbb{E}\exp\left\{-\sum_{l=1}^k\alpha_lE_v\left(\sum_{i=1}^lx_i\right)\right\}&=\widetilde{v}\left\{\lambda\left[1-\widetilde{N}\left(w(x,\alpha)\right)\right]\right\}\\&\cdot\prod_{l=1}^{k}\int_0^1\exp\left\{\lambda x_l\left[\widetilde{N}\left(s_l(x,\alpha,u)\right)-1\right]\right\}\rd u\,.\nonumber
		\end{align}
	\end{corollary}

	\subsection{Comparison with existing literature}
	Haviv and Ritov \cite{Haviv1998} considered the special case when $v$ is distributed according to the stationary distribution of the corresponding M/G/1  queue with an arrival rate $\lambda$ and a service distribution $B(\cdot)$. In this case, the expected value of $v$ is given by  
		\begin{equation}
		\mathbb{E}v=\frac{\lambda\mu_2}{2(1-\rho)}\,.
		\end{equation}
		Observe that an insertion of these formula into \eqref{eq: expectation1} provides exactly the same expression as in  \cite[Eqn.\ (7)]{Haviv1998}. Thus, in that sense, the formulae in this section may be considered as a natural generalization of this theorem, as in our framework $v$ can have any distribution. Importantly, the proof in the current work stems from other considerations than those which appeared in the original proof of  \cite{Haviv1998}. Moreover, Corollary~\ref{cor: LST} might be applied for the special case of $v$ which is distributed according to the stationary distribution of the corresponding M/G/1  system. This is a systematic approach to compute all externality moments in the model of  \cite{Haviv1998}.

	\section{Discussion and open problems}\label{sec: conclusion}
	The main contributions of this work lie in the introduction of the notion of the {externalities process} and in the derivation of various of its properties in the  case of a FCFS M/G/1  queue. The rest of this section includes a set of open problems, related to the research presented in this paper.
	
	\begin{enumerate}
		\item The current analysis is sensitive to the service discipline, in that it is FCFS-specific. Thus, it might be interesting to analyze the externalities process which corresponds to other service disciplines (e.g., preemptive ones) and examine the differences with respect to the results of the present paper. A particularly intriguing question concerns the characterization of the set of service disciplines for which the externalities process is convex.
		
		\item One could think about the externalities processes of more complex queues, e.g., G/G/1, M$_t$/G/1, etc. It is anticipated that in such cases the analysis is considerably more involved.       
		
		\item Consider the following natural multi-server version of the externalities process. Assume that there are $k$ servers and a Poisson arrival process of customers, where the service demands of the customers constitute a sequence of iid $k$-dimensional nonnegative random vectors which are independent from the arrival process. This defines $k$ coupled FCFS M/G/1  queues. Then, define a $k$-dimensional process such that its $i$-th ($1\leq i\leq k$) coordinate is the externalities process which is associated with the $i$-th queue. Also in this setup one would like to describe  the externalities process. A specific natural question is: Are there non-trivial assumptions on the $k$-dimensional service distributions under which we get an asymptotic independence of the externalities processes? 
		
		\item The L\'evy-driven queue, as analyzed in \cite{Debicki2015}, forms a class of storage models which can be seen as a natural generalization of the classic FCFS M/G/1  queue. A first question is: how should the externalities process be defined for such L\'evy queues? In particular, it is interesting to analyze whether there is a definition for which the results of the current work may be generalized relying on the machinery developed for L\'evy processes.   
	\end{enumerate}
	 
	 \section{Proofs}\label{sec: proofs}
	 \subsection{Proofs for Section \ref{sec: busy period}}

	 	\subsubsection*{Proof of Proposition \ref{prop: recursion}}
	 		Define a function 
	 		\begin{equation}\label{eq:f definition}
	 			f(x)=b\left[\lambda(1-x)\right])\ \ , \ \ x\in[0,1]
	 		\end{equation}
	 		and observe that for each $l\geq1$, the chain rule implies that
	 		\begin{equation}
	 			f^{(l)}(x)=(-\lambda)^lb^{(l)}\left[\lambda(1-x)\right]\,.
	 		\end{equation}
	 		In particular, when $x=1$, we get that $f^{(l)}(1)=\lambda^l\mu_l$. As a result, according to the Fa\'a di Bruno's formula, for each $k\geq1$ 
	 		\begin{align}
	 			\frac{\rd ^k}{\rd \alpha^k}b&\left\{\lambda\left[1-\widetilde{N}(\alpha)\right]\right\}\bigg|_{\alpha=0}=\frac{\rd ^k}{\rd \alpha^k}f\circ\widetilde{N}(\alpha)\bigg|_{\alpha=0}\\&=\sum_{m=1}^kf^{(m)}\left[\widetilde{N}(0)\right]\mathcal{B}_{k,m}\left[\widetilde{N}^{(1)}(0),\widetilde{N}^{(2)}(0),\ldots,\widetilde{N}^{(k-m+1)}(0)\right]\nonumber\\&=\sum_{m=1}^k\lambda^m\mu_m\check{\mathcal{B}}_{k,m}\,.\nonumber
	 		\end{align}
	 		Thus, observe that differentiating $n$ times (at zero) both sides of \eqref{eq: LST} with the general Leibniz rule yields that
	 		\begin{align}\label{eq: eta}
	 			(-1)^n\eta_n&=\sum_{k=0}^{n}\binom{n}{k}\left(\frac{\rd^{n-k}e^{-\alpha}}{\rd\alpha^{n-k}}\right)_{\alpha=0}\left(\frac{\rd ^k}{\rd \alpha^k}b\left\{\lambda\left[1-\widetilde{N}(\alpha)\right]\right\}\right)_{\alpha=0}\\&=(-1)^n+\sum_{k=1}^n\binom{n}{k}(-1)^{n-k}\sum_{m=1}^k\lambda^m\mu_m\check{\mathcal{B}}_{k,m}\,.\nonumber
	 		\end{align}
	 		Notice that $\eta_n$ appears in the right-hand side only in the term
	 		\begin{equation}
	 			\lambda\mu_1\check{\mathcal{B}}_{n,1}=\lambda\mu_1(-1)^n\eta_n=\rho(-1)^n\eta_n\,.
	 		\end{equation}
	 		This immediately yields the required recursive formula. $\blacksquare$
	 		
	 		\subsubsection*{Proof of Corollary \ref{cor: N moments}}
	 		Inserting $n=1$ into \eqref{eq: moment} immediately yields that $\eta_1=1/(1-\rho)$. In addition,
	 		\begin{align}
	 		    &\mathcal{B}_{1,1}(-\eta_1)=-\eta_1\,,\\&\mathcal{B}_{2,1}(-\eta_1,\eta_2)=\eta_2\,,\nonumber\\&\mathcal{B}_{2,2}(-\eta_1)=\eta_1^2\,,\nonumber\\&\mathcal{B}_{3,2}(-\eta_1,\eta_2)=-3\eta_1\eta_2\,,\nonumber\\&\mathcal{B}_{3,3}(-\eta_1)=-\eta_1^ 3\,.\nonumber
	 		\end{align}
	 		Thus, according to \eqref{eq: moment},
	 		\begin{align}
	 	     \eta_2&=\frac{1}{1-\rho}\cdot\left[1-\lambda\mu_1\mathcal{B}_{1,1}(-\eta_1)+\lambda^2\mu_2\mathcal{B}_{2,2}(-\eta_1)\right]\\&=\frac{1}{1-\rho}\cdot\left[1+\frac{2\rho}{1-\rho}+\frac{\lambda^2\mu_2}{(1-\rho)^2}\right]\,.\nonumber
	 		\end{align}
	 		In a similar fashion, we get that
	 		\begin{align}
	 		\eta_3&=-\frac{1}{1-\rho}\cdot\bigg\{-1-3\lambda\mu_1\mathcal{B}_{1,1}(-\eta_1)\\&\hspace{16mm}-3\left[\lambda\mu_1\mathcal{B}_{2,1}(-\eta_1,\eta_2)+\lambda^2\mu_2\mathcal{B}_{2,2}(-\eta_1)\right]\nonumber\\&\hspace{16mm}+\lambda^2\mu_2\mathcal{B}_{3,2}(-\eta_1,\eta_2)+\lambda^3\mu_3\mathcal{B}_{3,3}(-\eta_1)\bigg\}\nonumber\\&=\frac{1}{1-\rho}\left\{1+\frac{3\rho}{1-\rho}+3\left[\rho\eta_2+\frac{\lambda^2\mu_2}{(1-\rho)^2}\right]+\frac{3\lambda^2\mu_2\eta_2}{1-\rho}+\frac{\lambda^3\mu_3}{(1-\rho)^3}\right\}.\nonumber
	 		\end{align}
	 		
	 		\subsubsection*{Proof of Proposition \ref{prop: PGF}}
	 		$N(1)=b(\lambda)$ follows by differentiating both sides of \eqref{eq: PGF} at zero. Now, consider some $s\geq2$, then the general Leibniz rule and the Fa\'a di Bruno's formula (recall $f(\cdot)$, defined in \eqref{eq:f definition}) yield that
	 			\begin{align}
	 				&\frac{{\rd}^s}{\rd z^s}\left\{zb\left[\lambda\left(1-\hat{N}(z)\right)\right]\right\}\bigg|_{s=0}=\frac{\rd^{s-1}}{\rd z^{s-1}}b\left[\lambda\left(1-\hat{N}(z)\right)\right]\bigg|_{z=0}\\&=\frac{\rd ^{s-1}}{\rd z^{s-1}}\left[f\circ\hat{N}(z)\right]\bigg|_{z=0}\nonumber\\&=\sum_{m=1}^{s-1}f^{(m)}\left[\hat{N}(0)\right]\mathcal{B}_{s-1,m}\left[\hat{N}^{(1)}(0),\hat{N}^{(2)}(0),\ldots,\hat{N}^{(k-m+1)}(0)\right]\nonumber\\&=\sum_{m=1}^{s-1}(-\lambda)^mb^{(m)}(\lambda)\bar{\mathcal{B}}_{s-1,m}\,.\ \ \blacksquare\nonumber
	 			\end{align}
	 		
	 \subsection{Proofs for Section \ref{sec: definition}}
	 The proofs of Corollary \ref{cor: convex} and Theorem \ref{thm: decomposition} follow directly from the material presented in Section \ref{sec: definition}. Thus, we are now providing only the proof of Theorem \ref{thm: marginal}. 
	\subsubsection*{Proof of Theorem \ref{thm: marginal}}
		Observe that by definition of $\xi(x)$, for each $m\geq1$ and $y>0$,  
		\begin{equation}
			m\leq \xi(y)\Leftrightarrow\ \sum_{k=1}^mI_k\leq y\,.
		\end{equation}
		As a result, for every $x>0$ we have that:
		\begin{align}
			\int_0^x\dot{E}_v(y)\rd y-xM&=\int_0^x\sum_{m=1}^{\infty}N_m\textbf{1}_{\{m\leq \xi(y)\}}\,\rd y\\&=\sum_{m=1}^\infty N_m\int_0^x\textbf{1}_{\{m\leq \xi(y)\}}\rd y\nonumber\\&=\sum_{m=1}^\infty N_m\int_0^\infty\textbf{1}_{\left\{\sum_{k=1}^mI_{k}\leq y<x\right\}}\rd y\nonumber\\&=\sum_{m=1}^\infty N_{m}\left(x-\sum_{k=1}^mI_k\right)^+\,.\nonumber
		\end{align}
		With this identity at our disposal, the required result is a consequence of \eqref{eq: key identity}. $\blacksquare$
	
	 \subsection{Proofs of Section \ref{sec: moments}}
		\subsubsection*{Proof of Corollary \ref{cor: LST}}
	 	With the notations that have been used in Theorem \ref{thm: decomposition}, observe that \eqref{eq: decomposition 1} can be rephrased as follows:
	 	\begin{align}
	 	E_v\left(\sum_{i=1}^kx_i\right)&\stackrel{d}{=}\sum_{m=1}^{M_1}N_{m,1}(x_1+x_2+\ldots+x_k)\\&+\sum_{m=1}^{M_2}N_{m,2}(x_1U_{m,2}+x_2+x_3+\ldots+x_k)\nonumber\\&+\sum_{m=1}^{M_3}N_{m,3}(x_2U_{m,3}+x_3+x_4+\ldots+x_k)\nonumber\\&+\ldots+\sum_{m=1}^{M_{k+1}}N_{m,k+1}U_{m,k+1}x_k\,.\nonumber
	 	\end{align}
	 	Thus, we obtain that
	 	\begin{align}
	 	\sum_{l=1}^k\alpha_l&E_v\left(\sum_{i=1}^lx_i\right)\stackrel{d}{=}\sum_{m=1}^{M_1}N_{m,1}\left[\alpha_1x_1+\alpha_2(x_1+x_2)+\ldots+\alpha_k\left(\sum_{i=1}^kx_i\right)\right]\nonumber\\&+\sum_{m=1}^{M_2}N_{m,2}\left[\alpha_1x_1 U_{m,2}+\alpha_2\left(x_1 U_{m,2}+x_2\right)+\ldots+\alpha_k\left(x_1 U_{m,2}+\sum_{i=2}^kx_i\right)\right]\nonumber\\&+\sum_{m=1}^{M_3}N_{m,3}\left[\alpha_2x_2 U_{m,3}+\alpha_3\left(x_2 U_{m,3}+x_3\right)+\ldots+\alpha_k\left(x_2 U_{m,3}+\sum_{i=3}^kx_i\right)\right]\nonumber\\&+\ldots+\sum_{m=1}^{M_{k+1}}N_{m,k+1}\alpha_kx_kU_{m,k+1}\,.
	 	\end{align}
	 	Note that given $\mathcal{U}$, the sequences $\left\{N_{m,l}:m\geq1\right\}$, $1\leq l\leq k+1$ are independent. As a consequence, the result follows by conditioning and de-conditioning on $\mathcal{U}$ with an application of the LST formula of a  compound Poisson distribution. $\blacksquare$
	 \subsection{Proofs of Section \ref{sec: heavy traffic}}\label{proofTh5}
	 Since the proof of Theorem \ref{thm: externalities - Gaussian approx} includes an application of Theorem \ref{thm: CPP-Gaussian approx}, we start by providing the proof of Theorem \ref{thm: CPP-Gaussian approx}.
	 \subsubsection*{Proof of Theorem \ref{thm: CPP-Gaussian approx}}
	 The following well-known bound is  useful in the proof of Theorem \ref{thm: CPP-Gaussian approx}:
	 \begin{equation}\label{eq: exponent bound}
	 	\left|e^{iy}-1-iy+\frac{y^ 2}{2}\right|\leq \frac{|y|^3}{6}\ \ , \ \ \forall y\in\mathbb{R}\,.
	 \end{equation}
	 With this bound in hands, we  prove convergence of the finite-dimensional distributions as stated in the next lemma. For the proof, it is convenient to denote
	 \begin{equation}
	 	\theta_n\equiv\int_{-\infty}^\infty t\rd F_n(t)\ \ , \ \ \varrho_n\equiv\lambda_n\theta_n\ \ , \ \ n\geq1\,.
	 \end{equation}  
	 \begin{lemma}\label{lemma: finite dimension}
	 	The conditions of Theorem \ref{thm: CPP-Gaussian approx} imply that for every $1\leq d<\infty$ and $0\leq x_1<x_2<\ldots<x_d<\infty$,
	 	\begin{equation}
	 		\frac{1}{\sqrt{\lambda_n\mu_{2,n}}}\left[J_n(x_1),J_n(x_2),\ldots,J_n(x_d)\right]\stackrel{d}{\rightarrow}\mathcal{N}(0,\Sigma_x)\ \ \text{as}\ \ n\to\infty
	 	\end{equation}
	 	where $\Sigma$ is a covariance matrix such that $\Sigma_{ij}\equiv x_i\wedge x_j$ for every $1\leq i,j\leq d$.
	 \end{lemma}
	 \textbf{Proof:} To begin with, consider the special case $d=1$ and assume that for each $n\geq1$, $W_n$ is a random variable which is distributed according to $F_n(\cdot)$. Fix some $x\geq0$ and for each $n\geq1$ denote
	 	\begin{equation}
	 		\phi_{F_n}(y)\equiv\mathbb{E}e^{iyW_n}=\int_{-\infty}^ \infty e^{iyt}\rd F_n(t)\ \ , \ \ y\in\mathbb{R}\,.
	 	\end{equation}
	 	In particular, \eqref{eq: exponent bound} implies that for every $y\in\mathbb{R}$,
	 	\begin{align}
	 		\left|\phi_{F_n}(y)-1-iy\theta_n+\frac{y^2\sigma_n}{2}\right|&=\left|\mathbb{E}\left(e^{iyW_n}-1-iyW_n+\frac{y^2W_n^2}{2}\right)\right|\\&\leq\mathbb{E}\left|e^{iyW_n}-1-iyW_n+\frac{y^2W_n^2}{2}\right|\nonumber\\&\leq \frac{|y|^3\nu_n}{6}\,.\nonumber
	 	\end{align} 
	 	Thus, for a fixed $y\in\mathbb{R}$ and every $n\geq1$ we have 
	 	\begin{align}\label{eq: upper bound}
	 		&\left|\lambda_nx_1\left[\phi_{F_n}\left(\frac{y}{\sqrt{\lambda_n\sigma_n}}\right)-1\right]-i\frac{y\varrho_nx_1}{\sqrt{\lambda_n\sigma_n}}-\left(-\frac{y^2x_1}{2}\right)\right|\\&=\lambda_nx_1\left|\phi_{F_n}\left(\frac{y}{\sqrt{\lambda_n\sigma_n}}\right)-1-\left(i\frac{y\theta_n}{\sqrt{\lambda_n\sigma_n}}-\frac{y^2}{2\lambda_n}\right)\right|\leq\frac{|y|^3x_1\nu_n}{6\sqrt{\lambda_n\sigma^3_n}}\nonumber
	 	\end{align}
	 	and Condition \textbf{(II)} implies that the upper bound in  \eqref{eq: upper bound} converges to zero as $n\to\infty$. As a result, deduce that for every $y\in\mathbb{R}$,
	 	\begin{align}\label{eq: LST continuity}
	 		\lim_{n\to\infty}\mathbb{E}e^{iy\frac{J_n(x)}{\sqrt{\lambda_n\sigma_n}}}&=\lim_{n\to\infty}\exp\left\{\lambda_nx_1\left[\phi_{F_n}\left(\frac{y}{\sqrt{\lambda_n\sigma_n}}\right)-1\right]-i\frac{y\varrho_nx_1}{\sqrt{\lambda_n\sigma_n}}\right\}\nonumber\\&=\exp\left\{-\frac{y^2x_1}{2}\right\}
	 	\end{align}
	 	and hence the result follows (for $d=1$) by Levy's continuity theorem. 
	 	
	 	The next stage is to extend this result for $d>1$. To this end, for each $n\geq1$ define $d$ iid stochastic processes \[J^{(1)}_n(\cdot),J_n^{(2)}(\cdot),\ldots,J_n^{(d)}(\cdot)\] which have the same distribution as $J_n(\cdot)$. Since $J_n(\cdot)$ has stationary independent increments, then for each $n\geq1$ we get 
	 	\begin{align}
	 		\frac{1}{\sqrt{\lambda_n\sigma_n}}\left[J_n(x_1),J_n(x_2),\ldots,J_n(x_d)\right]
	 	\end{align}
	 	is distributed like 
	 	\begin{equation}\label{eq: finite-dimensional}
	 		\left[\frac{J^{(1)}_n(x_1)}{\sqrt{\lambda_n\sigma_n}},\frac{J^ {(2)}_n(x_2-x_1)}{\sqrt{\lambda_n\sigma_n}},\ldots,\frac{J^{(d)}_n(x_d-x_{d-1})}{\sqrt{\lambda_n\sigma_n}}\right]\Gamma
	 	\end{equation}
	 	where $\Gamma=\left[\gamma_{ij}\right]$ such that $\gamma_{ij}=\textbf{1}_{\{i\leq j\}}$ for every $1\leq i,j\leq d$.  The vector in \eqref{eq: finite-dimensional} consists of $d$ independent coordinates. Thus, according to the special case $d=1$, we deduce that 
	 	\begin{equation}
	 		\left[\frac{J^{(1)}_n(x_1)}{\sqrt{\lambda_n\sigma_n}},\frac{J^ {(2)}_n(x_2-x_1)}{\sqrt{\lambda_n\sigma_n}},\ldots,\frac{J^{(d)}_n(x_d-x_{d-1})}{\sqrt{\lambda_n\sigma_n}}\right]
	 	\end{equation}
	 	converges in distribution to
	 	\begin{equation}
	 		\mathcal{N}_d\left[\textbf{0},\text{diag}\left(x_1,x_2-x_1,\ldots,x_d-x_{d-1}\right)\right]
	 	\end{equation}
	 	as $n\to\infty$. Finally, it is readily verified that
	 	\begin{equation}
	 		\Gamma'\text{diag}\left(x_1,x_2-x_1,\ldots,x_d-x_{d-1}\right)\Gamma=\Sigma
	 	\end{equation}
	 	from which the result follows. $\blacksquare$
	 \newline\newline
	 Now, we are ready to prove the next lemma which is about validity of a tightness condition.
	 
	 \begin{lemma}\label{lemma: tightness}
	 	For every $\delta>0$ there exist $n_0\geq1$ and $\alpha,\beta>0$ such that 
	 	
	 	\begin{equation}
	 		\mathbb{P}\left\{\left|\frac{J_n(t)-J_n(s)}{\sqrt{\lambda_n\sigma_n}}\right|\leq\delta\right\}\geq\beta\ \ , \ \ \forall n\geq n_0\text{ and }t,s\geq0 \text{ s.t. }|t-s|<\alpha\,.
	 	\end{equation}
	 \end{lemma}
	 \textbf{Proof:} Fix some $\delta>0$ and take some $0\leq s<t$. Notice that $J_n(\cdot)$ is a process with stationary increments and $J_n(t-s)$ has a continuous distribution function. Therefore,  Lemma \ref{lemma: finite dimension} yields that
	 	\begin{align}\label{eq: tightness}
	 		\mathbb{P}\left\{\left|\frac{J_n(t)-J_n(s)}{\sqrt{\lambda_n\sigma_n}}\right|\leq\delta\right\}&=\mathbb{P}\left\{\left|\frac{J_n(t-s)}{\sqrt{\lambda_n\sigma_n}}\right|\leq\delta\right\}\\&=\mathbb{P}\left\{\frac{J_n(t-s)}{\sqrt{\lambda_n\sigma_n}}\leq\delta\right\}-\mathbb{P}\left\{\frac{J_n(t-s)}{\sqrt{\lambda_n\sigma_n}}\leq-\delta\right\}\nonumber\\&\xrightarrow{n\to\infty}\mathbb{P}\left\{\mathcal{N}(0,t-s)\leq\delta\right\}-\mathbb{P}\left\{\mathcal{N}(0,t-s)\leq-\delta\right\}\nonumber\\&=\mathbb{P}\left\{\left|\mathcal{N}(0,t-s)\right|<\delta\right\}\,.\nonumber
	 	\end{align}
	 	Clearly, the probability in the right-hand side of \eqref{eq: tightness} can be made sufficiently close to one by taking $s$ and $t$ which are close enough to each other, and hence the result follows. $\blacksquare$ \newline\newline
	 Finally, notice that for each $n\geq1$, $J_n(\cdot)$ has independent increments. Thus, using Lemma \ref{lemma: finite dimension} and Lemma \ref{lemma: tightness}, for each $k=1,2,\ldots$,  \cite[Theorem V.19]{Pollard2012} gives the required convergence on $\mathcal{D}[0,k]$ equipped with the uniform metric. Finally, to complete the proof of Theorem \ref{thm: CPP-Gaussian approx} it now suffices to apply  \cite[Theorem V.23]{Pollard2012}. $\blacksquare$

    \subsubsection*{Proof of Theorem \ref{thm: externalities - Gaussian approx}}
    Let $k>0$ and observe that Condition \textbf{(i)}, Condition \textbf{(ii)} and Condition \textbf{(iii)} allow us to apply Theorem \ref{thm: CPP-Gaussian approx} with the sequence $\left\{S_n(\cdot):n\geq1\right\}$ in order to deduce that
	 	\begin{equation}
	 		\frac{S_n(\cdot)}{\sqrt{\lambda_n\eta_{2,n}}}\Rightarrow W(\cdot)\ \ \text{as}\ \ n\to\infty
	 	\end{equation}
	 	when the convergence is in $\mathcal{D}[0,k+v]$ equipped with the uniform metric. Since the limit process is concentrated on $C[0,k+v]$, according to the representation theorem \cite[Theorem VI.13]{Pollard2012}, there is a probability space with a random processes \[\tilde{W}(\cdot)\stackrel{d}{=}W,\:\:\:\:\tilde{S}_{n}(\cdot)\stackrel{d}{=}S_{n,i}(\cdot),\:\:\:\:n\geq1,\] such that
	 	\begin{equation*}
	 		\sup_{0\leq y\leq k+v}\left|\frac{\tilde{S}_{n}(y)}{\sqrt{\lambda_n\eta_{2,n}}}-\tilde{W}(y)\right|\xrightarrow{n\to\infty}0
	 	\end{equation*}  
	 	with probability one. In particular, notice that
   \begin{align}
       \sup_{0\leq x\leq k}\left|\int_0^x\frac{\tilde{S}_{n}(v+y)}{\sqrt{\lambda_n\eta_{2,n}}}{\rm d}y-\int_0^x\tilde{W}(v+y){\rm d}y\right|\leq k\sup_{0\leq y\leq k+v}\left|\frac{\tilde{S}_{n}(y)}{\sqrt{\lambda_n\eta_{2,n}}}-\tilde{W}(y)\right|\,.
   \end{align}
   Since the RHS converges to zero with probability one, deduce that the process \begin{equation}
   x\mapsto\int_0^x\frac{\tilde{S}_{n}(v+y)}{\sqrt{\lambda_n\eta_{2,n}}}{\rm d}y\end{equation} admits weak convergence in $\mathcal{D}[0,k]$ equipped with the uniform metric to the process $x\mapsto\int_0^x\tilde{W}(v+y){\rm d}y$.  Especially, since $k$ is an arbitrary positive number and the limiting process is concentrated in $C[0,\infty)$, then this convergence can be extended to $\mathcal{D}[0,\infty)$ via  \cite[Theorem V.23]{Pollard2012}. Thus, the claim of Theorem \ref{thm: externalities - Gaussian approx} follows. $\blacksquare$
   
\subsubsection*{Proof of Proposition \ref{prop: sufficient conditions}}
	 Inserting  the expressions which appear in the statement of Corollary~\ref{cor: N moments} yields that for each $n\geq1$, \begin{equation}\label{eq: condition11}\frac{\eta_{3,n}}{\sqrt{\lambda_n\eta^3_{2,n}}}=
      \frac{1+\frac{3\rho_n}{1-\rho_n}+3\left[\rho_n\,\eta_{2,n}+\frac{\lambda_n^2\mu_{2,n}}{(1-\rho_n)^2}\right]+\frac{3\lambda_n^2\mu_{2,n}\,\eta_{2,n}}{1-\rho_n}+\frac{\lambda_n^3\mu_{3,n}}{(1-\rho_n)^3}}{\sqrt{\lambda_n(1-\rho_n)^2\eta_{2,n}^3}}\,,
	\end{equation}
	with
	\begin{equation}
	    \eta_{2,n}=\frac{1}{1-\rho_n}\cdot\left[1+\frac{2\rho_n}{1-\rho_n}+\frac{\lambda_n^2\mu_{2,n}}{(1-\rho_n)^2}\right].
	\end{equation}
	In addition, Jensen's inequality yields that 
	\begin{equation}\label{eq: cond1}
	    \limsup_{n\to\infty}\lambda_n^3\mu_{3,n}<\infty
	\end{equation}
	implies that 
	\begin{equation}\label{eq: cond2}
	    \limsup_{n\to\infty}\lambda_n^2\mu_{2,n}<\infty\,.
	\end{equation}
	Therefore, the assumption
	\begin{equation}\label{eq: limsup}
	    \limsup_{n\to\infty}\rho_n<1\,,
	\end{equation}
	implies that the nominator of \eqref{eq: condition11} is $\mathcal{O}(1)$ as $n\to\infty$. In addition, under the assumption \eqref{eq: limsup}, the denominator of \eqref{eq: condition11} is bounded from below by $(1-\rho_n)\sqrt{\lambda_n}$ which tends to $\infty$ as $n\to\infty$. The proof of the first statement follows immediately from these results.
	
	Due to the first statement, in order to prove the second statement, it is enough to consider the case when $\rho_n\uparrow1$ as $n\to\infty$. Notice that under the assumption 
	\begin{equation}
	    \liminf_{n\to\infty}\lambda^2\mu_{2,n}>0\,,
	\end{equation}
	the denominator of \eqref{eq: condition11} is  \[\Omega\left(\frac{\sqrt{\lambda_n}}{(1-\rho_n)^{7/2}}\right)\] as $n\to\infty$. In addition, due to \eqref{eq: cond1} and \eqref{eq: cond2}, the nominator of \eqref{eq: condition11} is $\mathcal{O}\left((1-\rho_n)^{-4}\right)$ as $n\to\infty$. Combining these results with the assumption that $\lambda_n(1-\rho_n)\rightarrow\infty$ as $n\to\infty$ completes the proof. $\blacksquare$
    \newline\newline
    {\textbf{Acknowledgment:} The authors would like to thank Moshe Haviv for his comments on an earlier version of the current work.}

\end{document}